\documentclass{article}

\usepackage{amsmath}
\usepackage{amssymb}
\usepackage{amsfonts}

\begin{document}

\renewcommand{\theequation}{\arabic{section}.\arabic{equation}}

\begin{center}
{\Large Exact controllability in minimal time of the Navier-Stokes periodic
flow in a 2D-channel}

\medskip

{\Large Gabriela Marinoschi}

\medskip

\textquotedblleft Gheorghe Mihoc-Caius Iacob\textquotedblright\ Institute of
Mathematical Statistics and

Applied Mathematics of the Romanian Academy,

Calea 13 Septembrie 13, Bucharest, Romania

gabriela.marinoschi@acad.ro

\bigskip
\end{center}

\noindent Abstract. This work is concerned with the necessary conditions of
optimality for a minimal time control problem $(P)$ related to the
linearized Navier-Stokes periodic flow in a 2D-channel, subject to a
boundary input which acts on the transversal component of the velocity. The
objective in this problem is the reaching of the laminar regime in a minimum
time, as well as its preservation after this time. The determination of the
necessary conditions of optimality relies on the analysis of intermediate
minimal time control problems $(P_{k})$ for the Fourier modes $"k"$
associated to the Navier-Stokes equations and on the proof of the maximum
principle for them. Also it is found that one can construct, on the basis of
the optimal controllers of problems $(P_{k}),$ a small time called here 
\textit{quasi minimal} and a boundary controller which realizes the required
objective in $(P).$

\bigskip


\bigskip


\section{Introduction}

\setcounter{equation}{0}

In this paper, we focus on the determination of the necessary conditions of
optimality for the linearized Navier-Stokes periodic flow in a channel,
driven in minimal time towards a stationary laminar regime by a boundary
control acting upon the transversal flow velocity.

Controllability of Navier-Stokes flow in finite small-time flow gave rise in
the literature to a set of reference works. In the following, we briefly
review some titles in the literature devoted to this subject. A maximum
principle for the time optimal control of the 2D Navier-Stokes equations is
presented in \cite{VB-97}. For some aspects concerning the Navier-Stokes
controllability and stabilization we refer the reader to the monographs \cite%
{VB-book-stab}, \cite{VB-control} and to the papers \cite{vb-stab-NA} and 
\cite{vb-stab-NS-SCL}, investigating the stabilization of the Navier-Stokes
flow in a channel by controllers with a vertical velocity observation which
acts on the normal component of velocity, and by noise wall normal
controllers, respectively. For techniques referring to minimal time
controllability results we also mention the papers \cite{Wang-Wang-2003}, 
\cite{Wang-Wang-2007}, \cite{GM-ESAIM}. Moreover, we cite the more recent
monograph \cite{Wang-book} where a detailed investigation of time optimal
control problems is developed.

The small-time global exact null controllability problem for the
Navier-Stokes equation was suggested by J.-L. Lions in \cite{Lions-61}. The
control was a source term supported within a small subset of the domain,
which is similar to controlling only a part of the boundary, with the
Dirichlet boundary condition on the uncontrolled part of the boundary. An
exhaustive presentation of various controllability problems, including also
that of Navier-Stokes equations, is found in \cite{Coron-book}. Here, we
shall indicate only a few titles related especially to the controllability
by using a boundary control. In\ \cite{AF-OI-31}, a small-time global exact
null controllability is proved for a control supported on the whole
boundary, while the paper \cite{Fursikov-Imanuvilov-1} is devoted to the
proof of the local exact controllability of the 2D Navier-Stokes system in a
bounded domain in the case when the control function is concentrated on the
whole boundary or on some part of it. The small-time global exact null
controllability for Navier-Stokes under an irrotational flow boundary
condition on the uncontrolled boundaries in a 2D rectangular domain it is
proved in \cite{Chapouly-13}. The exact boundary controllability of the
Navier-Stokes system where the controls are supported in a given open subset
of the boundary is provided in \cite{Rodrigues}. More recently, the paper 
\cite{Coron-2017} focuses on the small-time controllability presenting a new
method, which takes into account the boundary layer for getting the control
determination.

Let us consider the fluid flow in a $2$-$D$ infinitely long channel,
governed by the incompressible Navier-Stokes equations:

\begin{equation}
\begin{tabular}{ll}
$u_{t}-\nu \Delta u+uu_{x}+vu_{y}=\theta _{x},$ & $v_{t}-\nu \Delta
v+uv_{x}+vv_{y}=\theta _{y},$ \\ 
$u_{x}+v_{y}=0,$ &  \\ 
$u(t,x+2\pi ,y)=u(t,x,y),$ & $v(t,x+2\pi ,y)=v(t,x,y),$ \\ 
$u(t,x,0)=u(t,x,L)=0,$ & $v(t,x,0)=0,\mbox{ }v(t,x,L)=0,$ \\ 
$u(0,x,y)=u_{0},\mbox{ }v(0,x,y)=v_{0},$ & for $t\in \mathbb{R}%
_{+}=(0,\infty ),$ $x\in \mathbb{R},$ $y\in (0,L).$%
\end{tabular}
\label{0}
\end{equation}
Here, $(u,v)$ is the fluid velocity, $\theta $ is the pressure, the
subscripts $t,x,y$ represent the partial derivatives with respect to these
variables.

We also consider the steady-state flow with zero vertical velocity, governed
by (\ref{0}). This flow velocity turns out to be of the form $(U(y),0),$
where $U(y)=-\frac{a}{2\nu }\left( \frac{y^{2}}{L}-y\right) ,$ $a\in \mathbb{%
R}_{+}$ (see e.g., \cite{Temam})).

The problem we are concerned with is to steer the flow (\ref{0}) to the
stationary regime $U(y)$, within a minimal finite time, by means of a
boundary control $w$ acting at $y=L$ upon the transversal velocity component 
$v,$ namely by $v(t,x,L)=w(t,x).$ More precisely, the objective is to
characterize the boundary control $w$ which could force the flow $(u,v)$
starting from $(u_{0},v_{0})\neq (U(y),0)$ to reach the laminar regime $%
(U(y),0)$ at a minimal time and, moreover, to preserve it at this value
after that time. We stress that we are concerned with the determination of
the necessary conditions of optimality and not with the proof of the
controllability result. However, we succeed to prove that this action can be
done within a \textit{quasi-minimal} time, provided by optimal minimal times
for the problems in modes of the Fourier transform of the Navier-Stokes
linearized system. Further, we shall describe in detail the arguments.

We shall study this problem for the linearized flow around the laminar
steady-state $(U(y),0)$. Also, since the flow is periodic along the
longitudinal axis, we shall consider it on a period $(0,2\pi ).$ Thus, we
linearize (\ref{0}) around $(U(y),0)$ relying on the change of function $%
u\rightarrow u-U$ and continue to keep the same notation $u$ for the
linearized longitudinal velocity. Then, the linearized controlled system
reads%
\begin{equation}
\begin{tabular}{ll}
$u_{t}-\nu \Delta u+Uu_{x}+U_{y}v=\theta _{x},$ & $v_{t}-\nu \Delta
v+Uv_{x}=\theta _{y},$ \\ 
$u_{x}+v_{y}=0,$ &  \\ 
$u(t,2\pi ,y)=u(t,0,y),$ & $v(t,2\pi ,y)=v(t,0,y),$ \\ 
$u(t,x,0)=u(t,x,L)=0,$ & $v(t,x,0)=0,\mbox{ }v(t,x,L)=w(t,x),$ \\ 
$u(0,x,y)=u_{0},\mbox{ }v(0,x,y)=v_{0},$ & for $t\in \mathbb{R}_{+},$ $x\in
(0,2\pi ),$ $y\in (0,L).$%
\end{tabular}
\label{1}
\end{equation}

We express the flow controllability in minimal time by the problem%

\noindent$(P)$ $\mbox{\ \ \ \ \ \ \ \ \ Minimize }\Big\{ J(T,w)=T;\mbox{ }T>0,\mbox{ }w\in
H^{1}(0,T;L^{2}(0,2\pi )),  \mbox{ }w(0,x)=0, $

\hspace*{25mm} 
$\begin{array}{l}
\displaystyle

\mbox{\ \ \ \ \ \ \ \ \ \ \ \ }\int_{0}^{T}\int_{0}^{2\pi }\left\vert w_{t}(t,x)\right\vert ^{2}dxdt\leq
 \rho ^{2},\\
 \mbox{\ \ \ \ \ \ \ \ \ \ \ \ \ \ \ \ \ \ \ \ \ }u(T,x,y)=0,\mbox{ }%
v(T,x,y)=0\mbox{ a.e. }(x,y)\in (0,2\pi )\times (0,L) \Big\} , 
 \end{array}$

\noindent subject to system (\ref{1}), with $u_{0}\neq 0,$ $v_{0}\neq 0.$ In addition,
by resetting $w$ after $t=T$ one ensures that the null regime is preserved,
as we shall see.

It is obvious that the minimal time should be positive. Indeed, by absurd,
if $T=0,$ we would have $0=v(T)=v(0)=v_{0},$ and similarly for $u,$ which
contradicts the hypothesis. The requirement $w(0,x)=0$ is done especially
for technical purposes, but it is also in agreement with the fact that at
the initial time the boundary condition at $y=L$ is no-slip. Finally, it is
clear that in problem $(P)$ it is important to find information about the
controller $w$ only on the interval $(0,T)$ within which the objective is
reached. On the interval $(T,\infty )$ the function $w$ can take a whatever
value, in particular $0,$ and this choice will have the effect of preserving
the zero value for $v$ after the time $T$. That is why the property of $w$
of belonging to $H^{1}$ is required only for $t\in (0,T).$

Our purpose is to find the necessary conditions of optimality for $(P).$ To
this end, the following controllability assumption, which will allow the
derivation of an observability result for the adjoint system, will be in
effect:

\medskip

\noindent $(H)$ For each $(t_{0},T),$ $0\leq t_{0}<T<\infty ,$ and each $%
(u^{0},v^{0})\in (L^{2}(0,2\pi ;L^{2}(0,L)))^{2}$,

with $\left\Vert v^{0}\right\Vert _{L^{2}(0,2\pi ;(H^{2}(0,L))^{\ast })}\leq
1,$

there exists $w\in H^{1}(0,T;L^{2}(0,2\pi ))$ and $\gamma _{(t_{0},T)}>0,$ $%
\gamma _{(\cdot ,T)}\in L^{2}(0,T)$ with the properties%
\begin{equation*}
w(\tau )=0\mbox{ for }0\leq \tau \leq t_{0}<T,\mbox{ }\left(
\int_{0}^{T}\int_{0}^{2\pi }\left\vert w_{\tau }(\tau ,x)\right\vert
^{2}dxd\tau \right) ^{1/2}\leq \gamma _{(t_{0},T)},
\end{equation*}%
such that $u^{t_{0},w}(T,x,y)=0,$ $v^{t_{0},w}(T,x,y)=0$ a.e. $(x,y)\in
(0,2\pi )\times (0,L).$

\medskip

\noindent Here, $(u^{t_{0},w},v^{t_{0},w})$ is the solution to (\ref{1})
starting from $(u^{0},v^{0})\neq (0,0)$ at time $t=t_{0},$ and controlled by 
$w$ and $(H^{2}(0,L))^{\ast }$ is the dual of the space $H^{2}(0,L).$ We
note that $\gamma _{(t_{0},T)}$ depends on $(u^{0},v^{0}),$ on the interval $%
(t_{0},T),$ it is bounded on $(0,T-\delta )$ for all $\delta >0$ and $\gamma
_{(t_{0},T)}\rightarrow \infty $ as $t_{0}\rightarrow T.$ In fact $\gamma
_{(t_{0},T)}$ represents the controllability cost, which should be larger if
the objective is expected to be reached in a smaller time, but its
singularity is assumed to be square integrable (with respect to $t_{0}\in
(0,T)).$

In what concerns the controllability hypothesis $(H)$, its verification is
beyond the purpose of this work.

Next, we are going to describe the organization of the paper. It is
convenient to reduce problem $(P)$ to minimization problems for the Fourier
coefficients of the velocity. To this end, we write 
\begin{equation}
f(t,x,y)=\sum\limits_{k\in \mathbb{Z},\mbox{ }k\neq 0}f_{k}(t,y)e^{ikx},%
\mbox{ }f_{k}=\overline{f}_{-k},\mbox{ }f_{0}=0,\mbox{ }  \label{0-3}
\end{equation}%
(which ensures that $f$ is real), where $i=\sqrt{-1}\in \mathbb{C},$ the set
of complex numbers and $\overline{f}$ is the complex conjugate. The notation 
$f$ stands for $u,v,\theta ,$ and $f_{k}$ stands for $u_{k},v_{k},\theta
_{k}.$ Obviously, $f_{k}(t,y)\in \mathbb{C}$. Similarly, 
\begin{eqnarray}
w(t,x) &=&\sum\limits_{k\in \mathbb{Z},\mbox{ }k\neq 0}w_{k}(t)e^{ikx},\mbox{
}w_{k}=\overline{w}_{-k}  \label{0-4} \\
u_{0}(x,y) &=&\sum\limits_{k\in \mathbb{Z},\mbox{ }k\neq 0}u_{k0}(y)e^{ikx},%
\mbox{ }v_{0}(x,y)=\sum\limits_{k\in \mathbb{Z},\mbox{ }k\neq
0}v_{k0}(y)e^{ikx},\mbox{ }u_{k0}=\overline{u}_{-k0},\mbox{ }v_{k0}=%
\overline{v}_{-k0}.  \notag
\end{eqnarray}%
Replacing in (\ref{1}) the functions by their Fourier series and identifying
the coefficients we obtain the system%
\begin{eqnarray}
(u_{k})_{t}-\nu u_{k}^{\prime \prime }+(\nu k^{2}+ikU)u_{k}+U^{\prime }v_{k}
&=&ik\theta _{k},\mbox{ }  \label{2} \\
(v_{k})_{t}-\nu v_{k}^{\prime \prime }+(\nu k^{2}+ikU)v_{k} &=&\theta
_{k}^{\prime },  \notag \\
iku_{k}+v_{k}^{\prime } &=&0,\mbox{ }  \notag \\
u_{k}(t,0)=u_{k}(t,L)=v_{k}(t,0)=0,\mbox{ }v_{k}(t,L) &=&w_{k}(t),  \notag \\
u_{k}(0,y)=u_{k0},\mbox{ }v_{k}(0,y) &=&v_{k0}.  \notag
\end{eqnarray}%
For simplicity, we denote by the superscripts $^{\prime },^{\prime \prime
},^{\prime \prime \prime },^{\mbox{iv}}$ the first four partial derivatives
with respect to $y$ of the functions $u_{k}$ and $v_{k}.$

By the Parseval identity, in particular for%
\begin{equation}
\sum\limits_{k\in \mathbb{Z},\mbox{ }k\neq 0}\int_{0}^{T}\left\vert
(w_{k})_{t}(t)\right\vert ^{2}dt=\frac{1}{2\pi }\int_{0}^{2\pi }\left(
\int_{0}^{T}\left\vert w_{t}(t,x)\right\vert ^{2}dt\right) dx,  \label{0-6}
\end{equation}%
we can consider that 
\begin{equation*}
\int_{0}^{T}\left\vert (w_{k})_{t}(t)\right\vert ^{2}dt\leq \rho _{k}^{2},%
\mbox{ where}\sum\limits_{k\in \mathbb{Z},\mbox{ }k\neq 0}\rho _{k}^{2}\leq
\rho ^{2}.
\end{equation*}%
Eliminating $\theta _{k}$ between the first equations in (\ref{2}) by using 
\begin{equation}
v_{k}^{\prime }=-iku_{k}  \label{3}
\end{equation}%
we obtain the following equation in $v_{k}:$%
\begin{eqnarray}
(k^{2}v_{k}-v_{k}^{\prime \prime })_{t}+\nu v_{k}^{\mbox{iv}}-(2\nu
k^{2}+ikU)v_{k}^{\prime \prime }+(\nu k^{4}+ik^{3}U+ikU^{\prime \prime
})v_{k} &=&0,\mbox{ }  \label{4} \\
v_{k}(t,0)=0,\mbox{ }v_{k}(t,L) &=&w_{k}(t),  \label{4-1} \\
v_{k}^{\prime }(t,0)=v_{k}^{\prime }(t,L) &=&0,  \label{4-2} \\
v_{k}(0,y) &=&v_{k0}(y),  \label{4-3}
\end{eqnarray}%
for $t\in \mathbb{R}_{+},$ $y\in (0,L).$

Consequently, for each $k\in \mathbb{Z},$ $k\neq 0$ we can consider the
minimization problem for the mode "$k":$

\noindent$(P_k)$ $\qquad
\mbox{Minimize }\Big\{ J_{k}(T,w)=T;\mbox{ }T>0,\mbox{ }w\in
H^{1}(0,T),\mbox{ } w(0)=0, $
\hspace*{25mm} 

$\mbox{ \ \  \ \ \ \ \ \ \ \ \ \ \ \ \ \ \ \ \ \ \ \ \ \ \ \ \ \ \ \ \  } \mbox{ }\int_{0}^{T}\left\vert
w_{t}(t)\right\vert ^{2}dt\leq \rho _{k}^{2},\mbox{
}v_{k}(T,y)=0\mbox{ a.e. }y\in (0,L)\Big\} , $\\
\noindent subject to (\ref{4})-(\ref{4-3}), with $v_{k0}\neq 0.$

The controllability hypothesis $(H)$ and the Parseval identity provide for
each mode $"k"$ the following consequence:

\medskip

\noindent $(H_{k})$ For each $(t_{0},T),$ $0\leq t_{0}<T,$ and each initial
datum $v^{0}\in L^{2}(0,L),$ $v^{0}\neq 0,$ $\left\Vert v^{0}\right\Vert
_{(H^{2}(0,L))^{\ast }}\leq 1,$ there exists $w\in H^{1}(0,T)$ and $\gamma
_{(t_{0},T)}\in L^{1}(0,T),$ satisfying 
\begin{equation*}
w(\tau )=0\mbox{ for }0\leq \tau \leq t_{0},\mbox{ }\left(
\int_{0}^{T}\left\vert w_{\tau }(\tau )\right\vert ^{2}d\tau \right)
^{1/2}\leq \gamma _{(t_{0},T)},
\end{equation*}
such that $v^{t_{0},w}(T,y)=0,$ where $v^{t_{0},w}$ is the solution to (\ref%
{4})-(\ref{4-3}) starting from $v^{0}$ at time $t=t_{0}.$

\medskip

Here, we included the constant $\frac{1}{2\pi }$ in (\ref{0-6}) in $\gamma
_{(t_{0},T)}.$

As announced previously, the main part of the paper is directed to the
determination of the necessary conditions of optimality for $(P)$, which
will be deduced from those found for $(P_{k}).$

Here there is the structure of the paper. For technical reasons, by using an
appropriate variable transformation, we shall study instead of $(P_{k})$,
another problem $(\widehat{P_{k}})$ set on a fixed time interval. In Section
2, we shall prove, in Theorem 2.2, the well-posedness of the transformed
state system, and conclude with the existence of a solution $(T_{k}^{\ast
},w_{k}^{\ast })$ to $(\widehat{P_{k}})$. A characterization of the
optimality conditions cannot be done directly for problem $(\widehat{P_{k}}%
), $ so that we have to resort to an approximating problem $%
(P_{k,\varepsilon }),$ indexed along a small positive parameter $\varepsilon
,$ and prove the existence of a solution $(T_{k,\varepsilon }^{\ast
},w_{k,\varepsilon }^{\ast })$, in Theorem 3.1. A convergence result,
formally expressed by $(P_{k,\varepsilon })\rightarrow (\widehat{P_{k}}),$
as $\varepsilon \rightarrow 0,$ will be proved in Theorem 3.2. The latter
actually shows that if we fix an optimal pair $(T_{k}^{\ast },w_{k}^{\ast })$
in $(\widehat{P_{k}})$ we can recover it as a limit of a sequence of
solutions $(T_{k,\varepsilon }^{\ast },w_{k,\varepsilon }^{\ast })$ to $%
(P_{k,\varepsilon }).$ Based on these results, we proceed to the calculation
of the necessary conditions of optimality for the approximating minimization
problem $(P_{k,\varepsilon })$ in Proposition 3.4. They can be established
if $T_{k}^{\ast }$ is small enough and $\rho _{k}$ is chosen sufficiently
large. Appropriate fine estimates following by an observability result allow
to pass to the limit in the approximating optimality conditions to get those
corresponding to problem $(\widehat{P_{k}})$ in Theorem 4.1

Finally, by relying on the Fourier characterization of $u,$ $v$ and $w$ we
prove in Theorem 5.2 that, if $(P)$ has an admissible pair $(T_{\ast
},w_{\ast }),$ there exists $(T^{\ast },w^{\ast })$ which steers $%
(u_{0},v_{0})$ into $u(T^{\ast })=v(T^{\ast })=0$ and $T^{\ast }\leq T_{\ast
}.$ This pair is constructed on the basis of $(T_{k}^{\ast },w_{k}^{\ast })$
with $T_{k}^{\ast }$ minimal in problems $(P_{k}).$ For this reason we call
it a \textit{quasi minimal} time for $(P).$ However, it is not clear if it
is precisely the minimal one.

\medskip

\noindent \textbf{Notation.} Let $X_{\mathbb{R}}$ be a real Banach space and
let $T>0$. We denote by $X$ the complexified space $X_{\mathbb{R}}+iX_{%
\mathbb{R}}$ and by $L^{p}(0,T;X),$ $W^{1,p}(0,T;X),$ $C^{l}([0,T];X)$ the
complexified spaces containing functions of the form $f_{1}+if_{2},$ with $%
f_{1},f_{2}\in L^{2}(0,T;X_{\mathbb{R}}),$ $W^{k,p}(0,T;X_{\mathbb{R}})$ and 
$C^{l}([0,T];X_{\mathbb{R}}),$ respectively, for $p\in \lbrack 1,\infty ],$ $%
l\in \mathbb{N}.$

\noindent The space $W^{k,p}(0,T;X_{\mathbb{R}})=\{f\in L^{p}(0,T;X_{\mathbb{%
R}});$ $\frac{\partial ^{m}f}{\partial t^{m}}\in L^{p}(0,T;X_{\mathbb{R}}),$ 
$m=1,...,k\}.$

We shall use the standard Sobolev spaces $(H^{l}(0,1))_{\mathbb{R}},$ $%
(H_{0}^{1}(0,1))_{\mathbb{R}}$ and denote%
\begin{equation*}
H_{\mathbb{R}}:=(L^{2}(0,L))_{\mathbb{R}},\mbox{ }(H_{0}^{2}(0,L))_{\mathbb{R%
}}=\{f\in (H^{2}(0,L))_{\mathbb{R}};\mbox{ }f(0)=f(L)=f^{\prime
}(0)=f^{\prime }(L)=0\}.
\end{equation*}%
We have $(H_{0}^{2}(0,L))_{\mathbb{R}}\subset (H_{0}^{1}(0,L))_{\mathbb{R}%
}\subset H_{\mathbb{R}}\subset (H_{0}^{1}(0,L))_{\mathbb{R}}^{\ast }\subset
((H_{0}^{2}(0,L))_{\mathbb{R}})^{\ast }$ with compact injections, where $%
(H_{0}^{1}(0,L))_{\mathbb{R}}^{\ast }$ and $((H_{0}^{2}(0,L))_{\mathbb{R}%
})^{\ast }$ are the duals of $(H_{0}^{1}(0,L))_{\mathbb{R}}$ and $%
(H_{0}^{2}(0,L))_{\mathbb{R}},$ respectively. Their corresponding
complexified spaces $H_{0}^{2}(0,L),$ $H_{0}^{1}(0,L),$ $H,$ $%
(H_{0}^{1}(0,L))^{\ast },$ $(H_{0}^{2}(0,L))^{\ast }$ are defined as before,
and satisfy 
\begin{equation*}
H_{0}^{2}(0,L)\subset H_{0}^{1}(0,L)\subset H\subset (H_{0}^{1}(0,L))^{\ast
}\subset (H_{0}^{2}(0,L))^{\ast }
\end{equation*}%
with compact injections. Also, we define the spaces 
\begin{eqnarray}
V_{T} &:&=\{f\in H^{1}(0,T);\mbox{ }f(0)=0\},\mbox{ with the norm }%
\left\Vert w\right\Vert _{V_{T}}^{2}=\int_{0}^{T}\left\vert \dot{w}%
(t)\right\vert ^{2}dt,  \label{VT} \\
V_{1} &:&=\{f\in H^{1}(0,1);\mbox{ }f(0)=0\},\mbox{ with the norm }%
\left\Vert w\right\Vert _{V_{1}}^{2}=\int_{0}^{1}\left\vert \dot{w}%
(t)\right\vert ^{2}dt.  \notag
\end{eqnarray}%
We denote by $\left\vert \zeta \right\vert $ the norm of $\zeta =\zeta
_{1}+i\zeta _{2}\in \mathbb{C}$, the space of complex numbers. The scalar
product on $\mathbb{C}$ is defined as 
\begin{equation}
(a,b)_{\mathbb{C}}=a\overline{b},\mbox{ for }a,b\in \mathbb{C},\mbox{ with }%
\overline{b}\mbox{ the complex conjugate.}  \label{0-5}
\end{equation}%
The scalar product and norm in $H$ are defined by 
\begin{equation}
\left( \zeta ,z\right) _{H}=\int_{0}^{L}\zeta (y)\overline{z}(y)dy,\mbox{ }%
\left\Vert z\right\Vert _{H}=\left( \int_{0}^{L}\left\vert z(y)\right\vert
^{2}dy\right) ^{1/2},\mbox{ }\zeta ,z\in H.  \label{5}
\end{equation}%
The first and second derivatives of a function $w$ depending only on $t$
will be denoted by $\dot{w}$ and $\ddot{w}.$

\section{Problems $(P_{k})$ and $(\protect\widehat{P_{k}})$}

\setcounter{equation}{0}

In order to handle in a more convenient way the arguments in the proofs of
the next results, and especially for calculating the approximating necessary
conditions of optimality, we shall use a state system and a new minimization
problem for the modes $k$, on a fixed time interval, by making an
appropriate transformation in order to bring the interval $(0,T)$ into $%
(0,1).$ To this end, we set in the state system%
\begin{equation}
t=\widehat{t}T,\mbox{ }v_{k}(t,y)=\widehat{v_{k}}(\widehat{t}T,y):=\widehat{%
v_{k}}(\widehat{t},y),\mbox{ }\widehat{w_{k}}(t)=w_{k}(\widehat{t}),\mbox{ }
\label{t-cap}
\end{equation}%
such that $\widehat{t}\in \lbrack 0,1]$ when $t\in \lbrack 0,T].$

Then, the restriction $\left\Vert w\right\Vert _{V_{T}}\leq \rho _{k}$
becomes $\left\Vert \widehat{w}\right\Vert _{V_{1}}\leq \rho _{k}\sqrt{T},$
with $V_{1}$ defined in (\ref{VT}).

The state system (\ref{4})-(\ref{4-3}) is transformed into the appropriate
system for $\widehat{v_{k}}$ 
\begin{eqnarray}
(k^{2}\widehat{v_{k}}-\widehat{v_{k}}^{\prime \prime })_{\widehat{t}%
}+T\left( \nu \widehat{v_{k}}^{\mbox{iv}}-(2\nu k^{2}+ikU)\widehat{v_{k}}%
^{\prime \prime }+(\nu k^{4}+ik^{3}U+ikU^{\prime \prime })\widehat{v_{k}}%
\right) &=&0,\mbox{ }  \label{4'} \\
\widehat{v_{k}}(\widehat{t},0)=0,\mbox{ }\widehat{v_{k}}(\widehat{t},L) &=&%
\widehat{w_{k}}(\widehat{t}),  \label{4''} \\
\widehat{v_{k}}^{\prime }(\widehat{t},0)=\widehat{v_{k}}^{\prime }(\widehat{t%
},L) &=&0,  \label{4'''} \\
\widehat{v_{k}}(0,y) &=&v_{k0}(y),  \label{4''''}
\end{eqnarray}%
for $\widehat{t}\in (0,1),$ $y\in (0,L).$

In this way, problem $(P_{k})$ becomes $(\widehat{P_{k}})$ below:%

\noindent$(\widehat{P_{k}})$ $\quad\mbox{Minimize }\Big\{ J_{k}(T,\widehat w)=T;\ T>0,\ \widehat w\in V_{1},  \displaystyle\int_{0}^{1}\!\!\left\vert \dot{\widehat{w}} (t)\right\vert ^{2}dt\leq \rho _{k}^{2}T, \mbox{ } \widehat{v_{k}}(1,y)=0,\text{ a.e. }y\in (0,L)\Big\}$

\hspace*{75mm}
\medskip

\noindent subject to (\ref{4'})-(\ref{4''''}), with $v_{k0}\neq 0,$ where $V_{1}$ was
defined in (\ref{VT}).

The controllability hypothesis $(H_{k})$ will be correspondingly written on
the interval $(\widehat{t_{0}},1)$ and denoted by $(\widehat{H_{k}}).$

\medskip

\textbf{Note.} \textit{However, for not overloading the notation, we shall
skip in sections 2-4 the notation with the decoration }$"\symbol{94}"$ 
\textit{and will resume it in Theorem 5.2. Thus, in system} (\ref{4'})-(\ref%
{4''''}) \textit{we shall} \textit{write} $t,$\textit{\ }$v_{k},$\textit{\ }$%
w_{k}$\textit{, instead of }$\widehat{t},\widehat{v_{k}},\widehat{w_{k}}.$

\medskip

In this section we prove the well-posedness for the state system derived
from (\ref{4'})-(\ref{4''''}) and the existence of a solution to $(\widehat{%
P_{k}})$, which obviously imply the same results for (\ref{4})-(\ref{4-3})
and $(P_{k}).$

\medskip

We begin with some definitions. For each $k\in \mathbb{Z}\backslash \{0\}$
let us define the operators 
\begin{equation}
E_{0k}:D(E_{0k})\subset H\rightarrow H,\mbox{ }D(E_{0k})=H^{2}(0,L)\cap
H_{0}^{1}(0,L),\mbox{ }E_{0k}z:=k^{2}z-z^{\prime \prime }  \label{10-3}
\end{equation}%
and%
\begin{eqnarray}
F_{0k} &:&D(F_{0k})\subset H\rightarrow H,\mbox{ }D(F_{0k})=H^{4}(0,L)\cap
H_{0}^{2}(0,L),  \label{10-4} \\
F_{0k}z &=&\nu z^{\mbox{iv}}-(2\nu k^{2}+ikU)z^{\prime \prime }+(\nu
k^{4}+ik^{3}U+ikU^{\prime \prime })z.  \notag
\end{eqnarray}%
Since $E_{0k}$ is $m$-accretive, coercive, hence invertible, with the
inverse continuous on $H,$ we can define the operator 
\begin{equation}
A_{k}:=F_{0k}E_{0k}^{-1},\mbox{ }A_{k}:D(A_{k})\subset H\rightarrow H,\mbox{ 
}D(A_{k})=\{v\in H;\mbox{ }E_{0k}^{-1}v\in D(F_{0k})\}.  \label{10-2-4}
\end{equation}%
We also observe that,%
\begin{equation}
v\in D(A_{k})\mbox{ iff }v=E_{0k}\varphi ,\mbox{ for }\varphi \in
H^{4}(0,L)\cap H_{0}^{2}(0,L).  \label{iff}
\end{equation}

By Lemma 1 in \cite{vb-stab-NS-SCL} we know that $A_{k}$ is closed and
densely defined in $H,$ and $-A_{k}$ generates a $C_{0}$-analytic semigroup
on $H,$ that is, its resolvent has the property 
\begin{equation}
\left\Vert (\sigma I+A_{k})^{-1}f\right\Vert _{H}\leq \frac{\left\Vert
f\right\Vert _{H}}{\left\vert \sigma \right\vert -\sigma _{0}},\mbox{ for
all }f\in H\mbox{ and }\left\vert \sigma \right\vert >\sigma _{0}.
\label{resolvent}
\end{equation}

\medskip

\noindent \textbf{Definition 2.1. }We call a solution to\textbf{\ }(\ref{4'}%
)-(\ref{4''''}) a function 
\begin{equation*}
v_{k}\in C([0,1];H^{2}(0,L))\cap W^{1,2}(0,1;H^{2}(0,L))\cap
L^{2}(0,1;H^{4}(0,L)),
\end{equation*}
which satisfies (\ref{4'})-(\ref{4''''}) for a.a. $t>0.$

\medskip

\noindent \textbf{Theorem 2.2.} \textit{Let }$T>0$ \textit{and} $v_{k0}\in
H^{4}(0,L)\cap H_{0}^{2}(0,L),$ $w_{k}\in V_{1}.$\textit{\ Then, problem} (%
\ref{4'})-(\ref{4''''}) \textit{has a unique solution }%
\begin{equation}
v_{k}\in C([0,1];H^{2}(0,L))\cap W^{1,2}(0,1;H^{2}(0,L))\cap
L^{2}(0,1;H^{4}(0,L)),  \label{vk}
\end{equation}%
\textit{which satisfies the estimate}%
\begin{eqnarray}
&&\left\Vert v_{k}(t)\right\Vert _{H^{2}(0,L)}^{2}+\int_{0}^{1}\left\Vert 
\frac{dv_{k}}{dt}(t)\right\Vert _{H^{2}(0,L)}^{2}dt+T\int_{0}^{1}\left\Vert
v_{k}(t)\right\Vert _{H^{4}(0,L)}^{2}dt  \label{vk-est1} \\
&\leq &C\left( \left\Vert v_{k0}\right\Vert _{H^{4}(0,L)\cap
H_{0}^{2}(0,L)}^{2}+T\int_{0}^{1}\left\vert w_{k}(t)\right\vert
^{2}dt+\int_{0}^{1}\left\vert \dot{w}_{k}(t)\right\vert ^{2}dt\right) ,\mbox{
\textit{for all} }t\geq 0.  \notag
\end{eqnarray}%
\textit{The solution is continuous with respect to the data, that is, two
solutions }$(v_{k}^{1},v_{k}^{2})$\textit{\ corresponding to the data }$%
(v_{k0}^{1},w_{k}^{1})$\textit{\ and }$(v_{k0}^{2},w_{k}^{2})$\textit{\
satisfy the estimate} 
\begin{eqnarray}
&&\left\Vert (v_{k}^{1}-v_{k}^{2})(t)\right\Vert
_{H^{2}(0,L)}^{2}+\int_{0}^{1}\left\Vert \frac{d(v_{k}^{1}-v_{k}^{2})}{dt}%
(t)\right\Vert _{H^{2}(0,L)}^{2}dt  \label{vk-est2} \\
&&+T\int_{0}^{1}\left\Vert (v_{k}^{1}-v_{k}^{2})(t)\right\Vert
_{H^{4}(0,L)}^{2}dt  \notag \\
&\leq &C\left( \left\Vert (v_{k0}^{1}-v_{k0}^{2})\right\Vert
_{H^{4}(0,L)\cap H_{0}^{2}(0,L)}^{2}+\int_{0}^{1}\left( T\left\vert
(w_{k}^{1}-w_{k}^{2})(t)\right\vert ^{2}+\left\vert (\dot{w}_{k}^{1}-\dot{w}%
_{k}^{2})(t)\right\vert ^{2}\right) dt\right) ,\mbox{ }  \notag
\end{eqnarray}%
\textit{for all} $t\geq 0.$

\medskip

\noindent \textbf{Proof.} We recall that $V_{1}:=\{f\in H^{1}(0,1);$ $%
f(0)=0\}$. Let us introduce a function transformation in order to homogenize
the boundary conditions, namely 
\begin{equation}
\widetilde{v}_{k}(t,y)=v_{k}(t,y)-\beta (y)w_{k}(t),  \label{v-beta}
\end{equation}%
where%
\begin{equation}
\beta (y)=-\frac{2}{L^{3}}y^{3}+\frac{3}{L^{2}}y^{2}.  \label{beta}
\end{equation}%
This transformation is chosen such that $\widetilde{v}_{k,\varepsilon }(t,0)=%
\widetilde{v}_{k,\varepsilon }(t,L)=\widetilde{v}_{k,\varepsilon }^{\prime
}(t,0)=\widetilde{v}_{k,\varepsilon }^{\prime }(t,L)=0$. Equation (\ref{4'})
is transformed into 
\begin{eqnarray}
&&(k^{2}\widetilde{v}_{k}-\widetilde{v}_{k}^{\prime \prime })_{t}+T\left(
\nu \widetilde{v}_{k}^{\mbox{iv}}-(2\nu k^{2}+ikU)\widetilde{v}_{k}^{\prime
\prime }+(\nu k^{4}+ik^{3}U+ikU^{\prime \prime })\widetilde{v}_{k}\right)
\label{10-2} \\
&=&a_{k}(y)w_{k}(t)+Tb_{k}(y)\dot{w}_{k}(t),  \notag
\end{eqnarray}%
where%
\begin{equation}
a_{k}=-\left( \nu \beta ^{\mbox{iv}}-(2\nu k^{2}+ikU)\beta ^{\prime \prime
}+(\nu k^{4}+ik^{3}U+ikU^{\prime \prime })\beta \right) ,\mbox{ }%
b_{k}=-(k^{2}\beta -\beta ^{\prime \prime })  \label{ak-bk}
\end{equation}%
and $a_{k},$ $b_{k}\in C^{\infty }(0,L).$ We denote $\widetilde{\widetilde{v}%
}_{k}(t)=k^{2}\widetilde{v}_{k}(t)-\widetilde{v}_{k}^{\prime \prime }(t),$
for $t\in (0,1),$ and note that since $\widetilde{v}_{k}$ vanishes at the
boundaries, we have in fact 
\begin{equation}
\widetilde{\widetilde{v}}_{k}(t)=E_{0k}\widetilde{v}_{k}(t),\mbox{ for }t\in
(0,1).  \label{vL0}
\end{equation}%
Since $\widetilde{\widetilde{v}}_{k}(0)=(k^{2}v_{k0}-v_{k0}^{\prime \prime
})-(k^{2}\beta -\beta ^{\prime \prime })w_{k}(0)$ and $w_{k}(0)=0,$ equation
(\ref{10-2}) can be written as the equivalent Cauchy problem 
\begin{eqnarray}
\frac{d\widetilde{\widetilde{v}}_{k}}{dt}(t)+TA_{k}\widetilde{\widetilde{v}}%
_{k}(t) &=&Ta_{k}w_{k}(t)+b_{k}\dot{w}_{k}(t),\mbox{ a.e. }t\in (0,1),
\label{10-2-1} \\
\widetilde{\widetilde{v}}_{k}(0) &=&E_{0k}v_{k0},  \notag
\end{eqnarray}%
where $A_{k}=F_{0k}E_{0k}^{-1},$ by (\ref{10-2-4}). We recall that $-A_{k}$
generates a $C_{0}$-analytic semigroup and so the solution to (\ref{10-2-1})
is given by 
\begin{equation*}
\widetilde{\widetilde{v}}_{k}(t)=v_{1}(t)+v_{2}(t)+v_{3}(t),\mbox{ }
\end{equation*}%
where 
\begin{equation*}
v_{1}(t)=e^{-tTA_{k}}\widetilde{\widetilde{v}}_{k}(0),\mbox{ }%
v_{2}(t)=T\int_{0}^{t}(e^{-(t-s)TA_{k}}a_{k})w_{k}(s)ds,\mbox{ }%
v_{3}(t)=\int_{0}^{t}(e^{-(t-s)TA_{k}}b_{k})\dot{w}_{k}(s)ds.
\end{equation*}%
Since $v_{k0}\in H^{4}(0,L)\cap H_{0}^{2}(0,L)$, it follows that $\widetilde{%
\widetilde{v}}_{k}(0)=E_{0k}v_{k0}\in D(A_{k}).$ Then, using the existence
theorems for the solutions to equations with a $C_{0}$-analytic semigroup
(see e.g. \cite{Pazy}, Theorem 3.5, p. 114 and \cite{vbp-2012}, Proposition
1.148, p. 60), it follows for the first term that%
\begin{equation*}
v_{1}\in C([0,1];D(A_{k}))\cap C^{1}([0,1];H),
\end{equation*}%
\begin{equation*}
\left\Vert v_{1}(t)\right\Vert _{D(A_{k})}+\left\Vert \frac{dv_{1}}{dt}%
(t)\right\Vert _{H}\leq C\left\Vert E_{0k}v_{k0}\right\Vert _{D(A_{k})}\leq
C\left\Vert v_{k0}\right\Vert _{H^{4}(0,L)\cap H_{0}^{2}(0,L)},\mbox{ for
all }t\in \lbrack 0,1].
\end{equation*}%
The second term can be viewed as the solution to the Cauchy problem 
\begin{eqnarray*}
\frac{dv_{2}}{dt}(t)+TA_{k}v_{2}(t) &=&Ta_{k}w_{k}(t),\mbox{ a.e. }t\in
(0,1), \\
v_{2}(0) &=&0,
\end{eqnarray*}%
where $w_{k}a_{k}\in L^{2}(0,1;H),$ whence $v_{2}\in C([0,1];H)\cap
W^{1,2}(0,1;H)\cap L^{2}(0,1;D(A_{k}))$ and 
\begin{equation*}
\left\Vert v_{2}(t)\right\Vert _{H}^{2}+\int_{0}^{1}\left\Vert \frac{dv_{2}}{%
dt}(t)\right\Vert _{H}^{2}dt+T\int_{0}^{1}\left\Vert
A_{k}v_{2}(t)\right\Vert _{H}^{2}dt\leq CT\int_{0}^{1}\left\vert
w_{k}(t)\right\vert ^{2}dt.
\end{equation*}%
Similarly, $v_{3}\in C([0,1];H)\cap W^{1,2}(0,1;H)\cap L^{2}(0,1;D(A_{k}))$
and 
\begin{equation*}
\left\Vert v_{3}(t)\right\Vert _{H}^{2}+\int_{0}^{1}\left\Vert \frac{dv_{3}}{%
dt}(t)\right\Vert _{H}^{2}dt+T\int_{0}^{1}\left\Vert
A_{k}v_{3}(t)\right\Vert _{H}^{2}dt\leq C\int_{0}^{1}\left\vert \dot{w}%
_{k}(t)\right\vert ^{2}dt.
\end{equation*}%
Gathering the results for $v_{1},v_{2},v_{3},$ we get%
\begin{equation*}
\widetilde{\widetilde{v}}_{k}\in C([0,1];H)\cap W^{1,2}(0,1;H)\cap
L^{2}(0,1;D(A_{k})),
\end{equation*}%
\begin{eqnarray*}
&&\left\Vert \widetilde{\widetilde{v}}_{k}(t)\right\Vert
_{H}^{2}+\int_{0}^{1}\left\Vert \frac{d\widetilde{\widetilde{v}}_{k}}{dt}%
(t)\right\Vert _{H}^{2}dt+T\int_{0}^{1}\left\Vert A_{k}\widetilde{\widetilde{%
v}}_{k}(t)\right\Vert _{H}^{2}dt \\
&\leq &C\left( \left\Vert v_{k0}\right\Vert _{H^{4}(0,L)\cap
H_{0}^{2}(0,L)}^{2}+T\int_{0}^{1}\left\vert w_{k}(t)\right\vert
^{2}dt+\int_{0}^{1}\left\vert \dot{w}_{k}(t)\right\vert ^{2}dt\right) .
\end{eqnarray*}%
By (\ref{vL0}) we obtain 
\begin{equation*}
\widetilde{v}_{k}\in C([0,1];H^{2}(0,L)\cap H_{0}^{1}(0,L))\cap
W^{1,2}(0,1;H^{2}(0,L)\cap H_{0}^{1}(0,L))\cap L^{2}(0,1;H^{4}(0,L)\cap
H_{0}^{1}(0,L)),
\end{equation*}%
the last space being derived by $E_{0k}^{-1}(D(A_{k}))\subset H^{4}(0,L)\cap
H_{0}^{1}(0,L).$ Finally, by (\ref{v-beta}) we obtain (\ref{vk}). Estimate (%
\ref{vk-est1}) follows by the estimate for $\widetilde{\widetilde{v}}_{k}$.
Since the equation is linear we also get (\ref{vk-est2}). This implies still
the uniqueness. The proof is ended.\hfill $\square $

\medskip

Even if we are only interested in obtaining the conditions of optimality, we
also provide later, for the reader convenience, the proof of the existence
of a solution to $(P_{k}).$

To this end, we define an admissible pair for $(P)$ a pair $(T_{\ast
},w_{\ast })$ with $w_{\ast }\in H^{1}(0,T;L^{2}(0,2\pi )),$ $w_{\ast
}(0,y)=0,$ $\int_{0}^{T}\int_{0}^{2\pi }\left\vert w_{\ast
t}(t,x)\right\vert ^{2}dxdt\leq \rho ^{2}$ and $u(T,x,y)=0,$ $v(T,x,y)=0$.
By the Parseval identity this implies that $(T_{\ast },w_{\ast k})$ is an
admissible pair for $(\widehat{P_{k}}),$ with $w_{\ast k}$ being the mode $k$
of $w_{\ast }$. This will be resumed in Theorem 5.2. We are not concerned
here with the proof of the existence of an admissible pair. Related results
can be found in the literature already cited referring to the
controllability in small-time.

\medskip

\noindent \textbf{Theorem 2.3.} \textit{Let }$v_{k0}\in H^{4}(0,L)\cap
H_{0}^{2}(0,L),$ $v_{k0}\neq 0.$ \textit{If }$(\widehat{P_{k}})$\textit{\
has an admissible pair}, \textit{then it\ has at least a solution }$%
(T_{k}^{\ast },w_{k}^{\ast })$ \textit{with the corresponding optimal state }%
$v_{k}^{\ast }$\textit{. Moreover, let us set }$\widetilde{w_{k}}%
(t):=w_{k}^{\ast }(t)$ \textit{for} $t\in \lbrack 0,1],$ \textit{and }$%
\widetilde{w_{k}}(t)=0$\textit{\ for }$t\in (1,\infty ).$\textit{\ Then,} $%
v_{k}^{\ast }(t)=0$ \textit{for} $t>1.$

\medskip

\noindent \textbf{Proof. }Following the assumption before, $(\widehat{P_{k}}%
) $ has an infimum denoted $T_{k}^{\ast }$ which is positive. We consider a
minimizing sequence $(T_{k}^{n},w_{k}^{n})$ such that $T_{k}^{n}>0,$ $%
\left\Vert w_{k}^{n}\right\Vert _{V_{1}}\leq \rho _{k}\sqrt{T_{k}^{n}}$,
with $v_{k}^{n}(1,y)=0$ and 
\begin{equation}
T_{k}^{\ast }\leq J(T_{k}^{n},w_{k}^{n})=T_{k}^{n}\leq T_{k}^{\ast }+\frac{1%
}{n},\mbox{ }n\geq 1.  \label{J}
\end{equation}%
This yields $T_{k}^{n}\rightarrow T_{k}^{\ast }$ as $n\rightarrow \infty .$
Also, there exists $w_{k}^{\ast }\in H^{1}(0,1)$ such that, on a
subsequence, $w_{k}^{n}\rightarrow w_{k}^{\ast }$ weakly in $H^{1}(0,1),$
strongly in $C([0,1])$ by Arzel\`{a} theorem, and $\left\Vert w_{k}^{\ast
}\right\Vert _{V_{1}}\leq \rho _{k}\sqrt{T_{k}^{\ast }}.$ Thus, $%
w_{k}^{n}(0)\rightarrow w_{k}^{\ast }(0)=0.$ The solution to (\ref{4'})-(\ref%
{4''''}) corresponding to $(T_{k}^{n},w_{k}^{n})$ is denoted $v_{k}^{n},$
has the properties (\ref{vk})-(\ref{vk-est2}) with $T=T_{k}^{n}$ and $%
v_{k}^{n}(1,y)=0$. Thus, by a simple calculation, handling the property (\ref%
{vk-est2}) we get%
\begin{equation*}
v_{k}^{n}\rightarrow v_{k}^{\ast }\mbox{ strongly in }C[0,1];H^{2}(0,L))\cap
W^{1,2}(0,1;H^{2}(0,L))\cap L^{2}(0,1;H^{4}(0,L)).
\end{equation*}%
We can pass to the limit in (\ref{4'})-(\ref{4''''}) written for $%
(T_{k}^{n},w_{k}^{n})$ to obtain that $v_{k}^{\ast }$ is the solution to (%
\ref{4'})-(\ref{4''''}) corresponding to $(T_{k}^{\ast },w_{k}^{\ast }).$
Moreover, since $v_{k}^{n}\rightarrow v_{k}^{\ast }$ strongly in $C([0,1];H)$
we also have $v_{k}^{\ast }(1,y)=0,$ so that $(T_{k}^{\ast },w_{k}^{\ast })$
is optimal in $(\widehat{P_{k}}).$

Next, we prove the last assertion in the statement of the theorem. If $%
w_{k}^{\ast }$ is extended by $0$ on $(1,\infty ),$ the system for the
variable $\chi _{k}$ starting from the initial datum at $t=1$ reads 
\begin{eqnarray*}
(k^{2}\chi _{k}-\chi _{k}^{\prime \prime })_{t}+\nu \chi _{k}^{\mbox{iv}%
}-(2\nu k^{2}+ikU)\chi _{k}^{\prime \prime }+(\nu k^{4}+ik^{3}U+ikU^{\prime
\prime })\chi _{k} &=&0, \\
\chi _{k}(t,0)=0,\mbox{ }\chi _{k}(t,L)=0,\mbox{ }\chi _{k}^{\prime
}(t,0)=\chi _{k}^{\prime }(t,L) &=&0, \\
\chi _{k}(1,y)=v_{k}^{\ast }(1,y) &=&0,
\end{eqnarray*}%
for $(t,y)\in (1,\infty )\times (0,L)$. Obviously, it has the unique
solution $0$, which extends the solution $v_{k}^{\ast }$ on $(0,1).$\hfill $%
\square $

\medskip

We note that in problem $(\widehat{P_{k}}),$ the optimal state satisfies (%
\ref{4'})-(\ref{4''''}) with $T=T_{k,}^{\ast }$ and $w_{k}^{\ast }.$

\section{The approximating problem $(P_{k,\protect\varepsilon })$}

\setcounter{equation}{0}

In this section, we introduce an approximating minimization problem $%
(P_{k,\varepsilon })$, prove the existence of a solution, its convergence to 
$(\widehat{P_{k}})$ and determine the approximating necessary conditions of
optimality.

\noindent Let $(T_{k}^{\ast },w_{k}^{\ast })$ be a solution to $(\widehat{P_{k}})$ and
let $\varepsilon >0.$ We introduce the following approximating problem:

\noindent$(P_{k,\varepsilon})$\qquad\qquad 
$\mbox{Minimize}\Big\{ J_{k,\varepsilon }(T,w)=T+\displaystyle\frac{1}{2\varepsilon }%
\left\Vert (\sigma I+A_{k})^{-1}v_{k}(1)\right\Vert _{H}^{2}$

$\hspace*{40mm}+\displaystyle\frac{1}{2}%
\int_{0}^{1}\left\vert \int_{0}^{t}(w_{k}-w_{k}^{\ast })(\tau )d\tau
\right\vert ^{2}dt;\ T>0, \mbox{ },  w\in V_1,\ \left\Vert w\right\Vert _{V_{1}}\leq \rho _{k}\sqrt{T} \Big\} ,$
\\

\noindent subject to the approximating system (\ref{4'})-(\ref{4''''}). We underline
that $v_{k}(1)$ is $v_{k}(1,y).$

\medskip

\noindent \textbf{Theorem 3.1.} \textit{Let }$v_{k0}\in H^{4}(0,L)\cap
H_{0}^{2}(0,L),$ $v_{k0}\neq 0.$\textit{\ Then, problem} $(P_{k,\varepsilon
})$ \textit{has at least a solution }$(T_{k,\varepsilon }^{\ast
},w_{k,\varepsilon }^{\ast })$ \textit{with the corresponding optimal state }%
$v_{k,\varepsilon }^{\ast }$\textit{. }

\medskip

\noindent \textbf{Proof. }For\textbf{\ }$(P_{k,\varepsilon })$ we see that
there exists at least an admissible pair, which is $(T_{k}^{\ast
},w_{k}^{\ast }),$ the optimal pair in $(\widehat{P_{k}}).$ Then, $%
J_{k,\varepsilon }(T,w)$ is positive and so there exists $d_{\varepsilon
}=\inf J_{k,\varepsilon }(T,w)$ and it is positive. Indeed, by absurd if $%
J_{k,\varepsilon }(T,w)=0,$ then each term, including $T,$ should be equal
with 0. This implies that in the second term of $J_{k,\varepsilon },$ $%
v_{k}(T=0,y)=0,$ which contradicts $v_{k0}\neq 0.$

\noindent We consider a minimizing sequence $(T_{k,\varepsilon
}^{n},w_{k,\varepsilon }^{n})$ with $T_{k,\varepsilon }^{n}>0,$ $\left\Vert
w_{k,\varepsilon }^{n}(t)\right\Vert _{V_{1}}\leq \rho _{k}\sqrt{%
T_{k,\varepsilon }^{n}},$ satisfying 
\begin{equation}
d_{\varepsilon }\leq J_{k,\varepsilon }(T_{k,\varepsilon
}^{n},w_{k,\varepsilon }^{n})\leq d_{\varepsilon }+\frac{1}{n},\mbox{ }n\geq
1.  \label{14}
\end{equation}%
Hence, there exists $T_{k,\varepsilon }^{\ast }>0$ such that $%
T_{k,\varepsilon }^{n}\rightarrow T_{k,\varepsilon }^{\ast },$ as $%
n\rightarrow \infty .$ On a subsequence, denoted still by $n,$ we have $%
w_{k,\varepsilon }^{n}\rightarrow w_{k,\varepsilon }^{\ast },$ weakly in $%
V_{1}$ and strongly in $C([0,1]),$ so that $w_{k,\varepsilon
}^{n}(0)\rightarrow w_{k,\varepsilon }^{\ast }(0)=0$ and $\left\Vert
w_{k,\varepsilon }^{\ast }\right\Vert _{V_{1}}\leq \rho _{k}\sqrt{%
T_{k,\varepsilon }^{\ast }}.$ By (\ref{14}) 
\begin{equation}
\int_{0}^{t}(w_{k,\varepsilon }^{n}-w_{k}^{\ast })(\tau )d\tau \rightarrow
\int_{0}^{t}(w_{k,\varepsilon }^{\ast }-w_{k}^{\ast })(\tau )d\tau ,\mbox{
uniformly for all }t\in \lbrack 0,1],  \label{14-1}
\end{equation}%
according to Arzel\`{a} theorem, because the sequence $\left(
\int_{0}^{t}(w_{k,\varepsilon }^{n}-w_{k}^{\ast })(\tau )d\tau \right) _{n}$
is bounded in $L^{2}(0,1)$ and its derivative is bounded in $L^{2}(0,1),$
too.

Then, the state system (\ref{4'})-(\ref{4''''}) corresponding to $%
T_{k,\varepsilon }^{n}$ and $w_{k,\varepsilon }^{n}$ has, by Theorem 2.2, a
unique solution continuous in time on $[0,1].$ This solution $%
v_{k,\varepsilon }^{n}$ has the properties (\ref{vk})-(\ref{vk-est2}), with $%
T=T_{k,\varepsilon }^{n}$ and $w_{k,\varepsilon }^{n}.$ By (\ref{vk-est2})
we deduce that $v_{k,\varepsilon }^{n}\rightarrow v_{k,\varepsilon }^{\ast }$
strongly in $C([0,1];H^{2}(0,L))\cap W^{1,2}(0,1;H^{2}(0,L))\cap
L^{2}(0,1;H^{4}(0,L))$ as $n\rightarrow \infty .$ Passing to the limit in (%
\ref{4'})-(\ref{4''''}) written for $(T_{k,\varepsilon
}^{n},w_{k,\varepsilon }^{n})$ we get that $v_{k,\varepsilon }^{\ast }$ is
the solution to (\ref{4'})-(\ref{4''''}) corresponding to $T_{k,\varepsilon
}^{\ast }$ and $w_{k,\varepsilon }^{\ast }.$ Moreover, $(\sigma
I+A_{k})^{-1}v_{k,\varepsilon }^{n}(1)\rightarrow (\sigma
I+A_{k})^{-1}v_{k,\varepsilon }^{\ast }(1)$ strongly in $H.$ Passing to the
limit in (\ref{14}), as $n\rightarrow \infty ,$ we get on the basis of the
previous convergences and of the weakly lower semicontinuity of the norms,
that $J_{k,\varepsilon }(T_{k,\varepsilon }^{\ast },w_{k,\varepsilon }^{\ast
})=d_{\varepsilon },$ that is $(T_{k,\varepsilon }^{\ast },w_{k,\varepsilon
}^{\ast })$ is an optimal controller in $(P_{k,\varepsilon }).$ \hfill $%
\square $

\medskip

We note that in problem $(P_{k,\varepsilon }),$ the optimal state satisfies (%
\ref{4'})-(\ref{4''''}) with $T=T_{k,\varepsilon }^{\ast }$ and $%
w_{k,\varepsilon }^{\ast }.$

\medskip

\noindent \textbf{Theorem 3.2.} \textit{Let }$(T_{k,\varepsilon }^{\ast
},w_{k,\varepsilon }^{\ast },v_{k,\varepsilon }^{\ast })$ \textit{be optimal
in }$(P_{k,\varepsilon })$\textit{\ and }$(T_{k}^{\ast },w_{k}^{\ast
},v_{k}^{\ast })$\textit{\ be optimal in }$(\widehat{P_{k}}).$ \textit{Then, 
}%
\begin{equation}
T_{k,\varepsilon }^{\ast }\rightarrow T_{k}^{\ast },\mbox{ }w_{k,\varepsilon
}^{\ast }\rightarrow w_{k}^{\ast }\mbox{ \textit{weakly in} }H^{1}(0,1)\mbox{
\textit{and strongly in} }C([0,1]),  \label{50}
\end{equation}%
\begin{equation}
v_{k,\varepsilon }^{\ast }\rightarrow v_{k}^{\ast }\mbox{ \textit{strongly in%
} }C([0,1];H^{2}(0,L))\cap W^{1,2}(0,1;H^{4}(0,L))\cap L^{2}(0,1;H^{4}(0,L)).
\label{51}
\end{equation}

\medskip

\noindent \textbf{Proof. }Let $(T_{k,\varepsilon }^{\ast },w_{k,\varepsilon
}^{\ast },v_{k,\varepsilon }^{\ast })$\textit{\ }be optimal in $%
(P_{k,\varepsilon })$ and denote 
\begin{equation*}
h_{k,\varepsilon }^{\ast }(t)=\int_{0}^{t}(w_{k,\varepsilon }^{\ast
}-w_{k}^{\ast })(\tau )d\tau .
\end{equation*}%
The fact that $(T_{k,\varepsilon }^{\ast },w_{k,\varepsilon }^{\ast })$ is
optimal in $(P_{k,\varepsilon })$ implies that 
\begin{eqnarray}
&&J_{k,\varepsilon }(T_{k,\varepsilon }^{\ast },w_{k,\varepsilon }^{\ast
})=T_{k,\varepsilon }^{\ast }+\frac{1}{2\varepsilon }\left\Vert (\sigma
I+A_{k})^{-1}v_{k,\varepsilon }^{\ast }(1)\right\Vert _{H}^{2}+\frac{1}{2}%
\int_{0}^{1}\left\vert h_{k,\varepsilon }^{\ast }(t)\right\vert ^{2}dt
\label{52-0} \\
&\leq &J_{k,\varepsilon }(T,w)=T+\frac{1}{2\varepsilon }\left\Vert (\sigma
I+A_{k})^{-1}v(1)\right\Vert _{H}^{2}+\frac{1}{2}\int_{0}^{1}\left\vert
\int_{0}^{t}(w-w_{k}^{\ast })(\tau )d\tau \right\vert ^{2}dt  \notag
\end{eqnarray}%
for any $T>0$ and $w\in V_{1},$ $\left\Vert w\right\Vert _{V_{1}}\leq \rho
_{k}\sqrt{T},$ where $v$ is the solution to the state system corresponding
to $(T,w).$ Let us set in (\ref{52-0}), $T=T_{k}^{\ast }$ and $w=w_{k}^{\ast
},$ the chosen optimal controller in $(\widehat{P_{k}}).$ Thus, the second
and the last terms on the right-hand side of (\ref{52-0}) vanish and 
\begin{equation}
J_{k,\varepsilon }(T_{k,\varepsilon }^{\ast },w_{k,\varepsilon }^{\ast
})=T_{k,\varepsilon }^{\ast }+\frac{1}{2\varepsilon }\left\Vert (\sigma
I+A_{k})^{-1}v_{k,\varepsilon }^{\ast }(1)\right\Vert _{H}^{2}+\frac{1}{2}%
\int_{0}^{1}\left\vert h_{k,\varepsilon }^{\ast }(t)\right\vert ^{2}dt\leq
T_{k}^{\ast }.  \label{52-1}
\end{equation}%
Then, $T_{k,\varepsilon }^{\ast }\rightarrow T_{k}^{\ast \ast },$ and on a
subsequence denoted still by $\varepsilon ,$ we have $w_{k,\varepsilon
}^{\ast }\rightarrow w_{k}^{\ast \ast }$ weakly in $V_{1},$ strongly in $%
C([0,1]),$ and $w_{k,\varepsilon }^{\ast }(0)\rightarrow w_{k}^{\ast \ast
}(0)=0.$ Also, $\left\Vert w_{k}^{\ast \ast }\right\Vert _{V_{1}}\leq \rho
_{k}\sqrt{T_{k}^{\ast \ast }}.$

The solution $v_{k,\varepsilon }^{\ast \ast }$ corresponding to $%
(T_{k}^{\ast \ast },w_{k,\varepsilon }^{\ast \ast })$ exists, it is unique,
according to Theorem 2.2 and has the properties (\ref{vk})-(\ref{vk-est2})
with $T=T_{k,\varepsilon }^{\ast }$ and $w_{k,\varepsilon }^{\ast }$.
Therefore, by handling some calculations based on (\ref{vk-est2}) we get 
\begin{equation}
v_{k,\varepsilon }^{\ast }\rightarrow v_{k}^{T_{k}^{\ast \ast },w_{k}^{\ast
\ast }}:=v_{k}^{\ast \ast }\mbox{ strongly in }W^{1,2}(0,1;H^{2}(0,L))\cap
L^{2}(0,1;H^{4}(0,L))  \label{52-2}
\end{equation}%
where $v_{k}^{\ast \ast }$ is the solution to (\ref{4'})-\ref{4''''})
corresponding to $(T_{k,\varepsilon }^{\ast \ast },w_{k,\varepsilon }^{\ast
\ast })$. By (\ref{52-1}) we have%
\begin{equation}
\left\Vert (\sigma I+A_{k})^{-1}v_{k,\varepsilon }^{\ast }(1)\right\Vert
_{H}^{2}\leq 2\varepsilon T_{k}^{\ast },  \label{e1}
\end{equation}%
so that 
\begin{equation}
\lim_{\varepsilon \rightarrow 0}\left\Vert (\sigma
I+A_{k})^{-1}v_{k,\varepsilon }^{\ast }(1)\right\Vert _{H}^{2}=0  \label{e2}
\end{equation}%
which implies that $v_{k,\varepsilon }^{\ast }(1,\cdot )\rightarrow 0$
strongly in $H.$ On the other hand, by (\ref{52-2}), $v_{k,\varepsilon
}^{\ast }(1)\rightarrow v_{k}^{\ast \ast }(1),$ so that $v_{k}^{\ast \ast
}(1)=0.$ Again by (\ref{52-1}), $T_{k,\varepsilon }^{\ast }\leq
J_{k,\varepsilon }(T_{\varepsilon }^{\ast },w_{\varepsilon }^{\ast })\leq
T_{k}^{\ast },$ implying at limit that $T_{k}^{\ast \ast }\leq T_{k}^{\ast
}. $ Since $T_{k}^{\ast \ast }$ and $w_{k}^{\ast \ast }$ satisfy the
restrictions required in problem $(\widehat{P_{k}}),$ that is $T_{k}^{\ast
\ast }>0,$ $\left\Vert w_{k}^{\ast \ast }\right\Vert _{V_{1}}\leq \rho _{k}%
\sqrt{T_{k}^{\ast \ast }},$ and $v_{k}^{\ast \ast }(1)=0,$ recalling that $%
T_{k}^{\ast }$ is the infimum in $(\widehat{P_{k}})$ it follows that $%
T_{k}^{\ast \ast }=T_{k}^{\ast }.$

Again by (\ref{52-1}) we see that $T_{k}^{\ast }+\limsup\limits_{\varepsilon
\rightarrow 0}\frac{1}{2\varepsilon }\left\Vert (\sigma
I+A_{k})^{-1}v_{k,\varepsilon }^{\ast }(1)\right\Vert _{H}^{2}\leq
T_{k}^{\ast },$ which implies 
\begin{equation}
\limsup_{\varepsilon \rightarrow 0}\frac{1}{2\varepsilon }\left\Vert (\sigma
I+A_{k})^{-1}v_{k,\varepsilon }^{\ast }(1)\right\Vert _{H}^{2}=0.
\label{52-4}
\end{equation}%
Also, $T_{k}^{\ast }+\limsup\limits_{\varepsilon \rightarrow 0}\frac{1}{2}%
\int_{0}^{1}\left\vert h_{k,\varepsilon }^{\ast }(t)\right\vert ^{2}dt\leq
T_{k}^{\ast },$ implying that 
\begin{equation}
\limsup_{\varepsilon \rightarrow 0}\int_{0}^{1}\left\vert h_{k,\varepsilon
}^{\ast }(t)\right\vert ^{2}dt=0.  \label{52-5}
\end{equation}%
Therefore, 
\begin{equation}
h_{k,\varepsilon }^{\ast }\rightarrow 0\mbox{ strongly in }L^{2}(0,1),\mbox{
as }\varepsilon \rightarrow 0  \label{52-6}
\end{equation}%
and so it follows that 
\begin{equation*}
h_{k,\varepsilon }^{\ast }(t)=\int_{0}^{t}(w_{k,\varepsilon }^{\ast
}-w_{k}^{\ast })(\tau )d\tau \rightarrow \int_{0}^{t}(w_{k}^{\ast \ast
}-w_{k}^{\ast })(\tau )d\tau =0,\mbox{ for all }t\in \lbrack 0,1],\mbox{ as }%
\varepsilon \rightarrow 0.
\end{equation*}%
Thus, we get $w_{k}^{\ast \ast }=w_{k}^{\ast }$ for all $[0,1]$ and so $%
v_{k}^{\ast \ast }(t)=v_{k}^{\ast \ast }(t)$ for all $t\in \lbrack 0,1].$

\noindent Based on the previous convergences, we pass to the limit in (\ref%
{52-0}) and conclude that $\lim\limits_{\varepsilon \rightarrow
0}J_{k,\varepsilon }(T_{k,\varepsilon }^{\ast \ast },w_{k,\varepsilon
}^{\ast \ast })=T_{k}^{\ast }=J_{k}(T_{k}^{\ast },w_{k}^{\ast })=T_{k}^{\ast
}$. \hfill $\square $

\subsection{Systems in variations and the adjoint system for $(P_{k,\protect%
\varepsilon })$}

Since the minimization problem depends on two controllers, $T_{k,\varepsilon
}^{\ast }$ and $w_{k,\varepsilon }^{\ast },$ we shall study separate
variations with respect to them. First, let us keep $T_{k,\varepsilon
}^{\ast }$ fixed and give variations to $w_{k,\varepsilon }^{\ast }.$ We
shall obtain a first system in variations.

\noindent Let $(T_{k,\varepsilon }^{\ast },w_{k,\varepsilon }^{\ast })$ be
an optimal controller in $(P_{k,\varepsilon }).$ For $\lambda >0,$ we set%
\begin{equation}
w_{k,\varepsilon }^{\lambda }=w_{k,\varepsilon }^{\ast }+\lambda \omega ,%
\mbox{ where }\omega =\widetilde{w}-w_{k,\varepsilon }^{\ast },\mbox{ }%
\widetilde{w}\in V_{1}.  \label{53}
\end{equation}%
We note that the state system satisfies (\ref{4'})-(\ref{4''''}) with $%
T=T_{k,\varepsilon }^{\ast }$ and $w_{k,\varepsilon }^{\ast }.$ We define $%
Y_{\lambda }=\frac{v_{k,\varepsilon }^{T_{\varepsilon }^{\ast
},w_{\varepsilon }^{\lambda }}-v_{k,\varepsilon }^{\ast }}{\lambda },$ where 
$v_{k,\varepsilon }^{T_{k,\varepsilon }^{\ast },w_{k,\varepsilon }^{\lambda
}}$ is the solution to (\ref{4'})-(\ref{4''''}) corresponding to $%
T_{k,\varepsilon }^{\ast }$ and $w_{k,\varepsilon }^{\lambda }.$ Taking into
account the estimates of Theorem 2.2, we deduce by a straightforward
calculation that 
\begin{equation*}
Y_{\lambda }\rightarrow Y\mbox{ as }\lambda \rightarrow 0,\mbox{ strongly in 
}C([0,1];H^{2}(0,L))\cap W^{1,2}(0,1;H^{2}(0,L))\cap L^{2}(0,1;H^{4}(0,L)),
\end{equation*}%
and that $Y$ satisfies the equations%
\begin{eqnarray}
(k^{2}Y-Y^{\prime \prime })_{t}+T_{k,\varepsilon }^{\ast }\left( \nu Y^{%
\mbox{iv}}-(2\nu k^{2}+ikU)Y^{\prime \prime }+(\nu k^{4}+ik^{3}U+ikU^{\prime
\prime })Y\right) &=&0,\mbox{ }  \label{54} \\
Y(t,0)=0,\mbox{ }Y(t,L)=\omega ,\mbox{ }Y^{\prime }(t,0)=Y^{\prime }(t,L)
&=&0,  \label{54-1} \\
Y(0,y) &=&0,  \label{54-3}
\end{eqnarray}%
for $(t,y)\in (0,1)\times (0,L).$ Moreover, following Theorem 2.2, we state
that (\ref{54})-(\ref{54-3}) has a unique solution 
\begin{equation}
Y\in C([0,1];H^{2}(0,L))\cap W^{1,2}(0,1;H^{2}(0,L))\cap
L^{2}(0,1;H^{4}(0,L)).  \label{55}
\end{equation}

We introduce the dual system%
\begin{equation}
-(k^{2}p_{k,\varepsilon }-p_{k,\varepsilon }^{\prime \prime
})_{t}+T_{k,\varepsilon }^{\ast }\left( \nu p_{k,\varepsilon }^{\mbox{iv}%
}-(2\nu k^{2}+ikU)p_{k,\varepsilon }^{\prime \prime }-2ikU^{\prime
}p_{k,\varepsilon }^{\prime }+(\nu k^{4}+ik^{3}U)p_{k,\varepsilon }\right)
=0,\mbox{ }  \label{56-1}
\end{equation}%
\begin{equation}
p_{k,\varepsilon }(t,0)=p_{k,\varepsilon }(t,L)=p_{k,\varepsilon }^{\prime
}(t,0)=p_{k,\varepsilon }^{\prime }(t,L)=0,  \label{56-2}
\end{equation}%
\begin{equation}
k^{2}p_{k,\varepsilon }(1,y)-p_{k,\varepsilon }^{\prime \prime
}(1,y)=(\sigma I+A_{k})^{-2}\left( \frac{1}{\varepsilon }\overline{%
v_{k,\varepsilon }^{\ast }}(1)\right) ,  \label{56-3}
\end{equation}%
for $(t,y)\in (0,1)\times (0,L).$

Now, in order to simplify the writing in the next calculations, we use a
formal notation 
\begin{equation}
Ev:=k^{2}v-v^{\prime \prime },\mbox{ }Fz:=\nu v^{\mbox{iv}}-(2\nu
k^{2}+ikU)v^{\prime \prime }+(\nu k^{4}+ik^{3}U+ikU^{\prime \prime })v
\label{L-V}
\end{equation}%
for $v\in H^{4}(0,L)$ and rewrite the state equation (\ref{4'})-(\ref{4''''}%
) for the solution $v_{k,\varepsilon }^{\ast }$ corresponding to $%
w_{k,\varepsilon }^{\ast }$ as 
\begin{equation}
E(v_{k,\varepsilon }^{\ast })_{t}(t)+T_{k,\varepsilon }^{\ast
}Fv_{k,\varepsilon }^{\ast }(t)=0,\mbox{ a.e. }t\in (0,1),  \label{75}
\end{equation}%
\begin{equation}
v_{k,\varepsilon }^{\ast }(0)=v_{k0}.  \label{76}
\end{equation}

\medskip

\noindent \textbf{Proposition 3.3. }\textit{The adjoint system} (\ref{56-1}%
)-(\ref{56-3}) \textit{has a unique solution}\textbf{\ }%
\begin{equation}
p_{k,\varepsilon }\in C([0,1];H^{6}(0,L)\cap H_{0}^{2}(0,L))\cap
C^{1}([0,1];H^{4}(0,L)\cap H_{0}^{2}(0,L)).  \label{57}
\end{equation}

\medskip

\noindent \textbf{Proof. }We recall\textbf{\ }the definition (\ref{10-3})
and introduce, similarly to (\ref{10-4}), the following operator:%
\begin{eqnarray}
F_{0k}^{\ast } &:&D(F_{0k}^{\ast })=D(F_{0k})\subset H\rightarrow H,
\label{57-2} \\
F_{0k}^{\ast }z &=&\nu z^{\mbox{iv}}-(2\nu k^{2}+ikU)z^{\prime \prime
}-2ikU^{\prime }z^{\prime }+(\nu k^{4}+ik^{3}U)z.  \notag
\end{eqnarray}%
We can interpret $F_{0k}^{\ast }$ as the dual of $F$ in the sense of
distributions, that is $(Fv,\overline{\varphi })=(\overline{v},F_{0k}^{\ast
}\varphi ),$ for $\varphi \in C_{0}^{\infty }(0,L)$ and $v\in H^{4}(0,L)$.
System (\ref{56-1})-(\ref{56-3})\textbf{\ }can be written%
\begin{equation}
-E_{0k}(p_{k,\varepsilon })_{t}(t)+T_{k,\varepsilon }^{\ast }F_{0k}^{\ast
}p_{k,\varepsilon }(t)=0,\mbox{ a.e. }t\in (0,1),  \label{57-0}
\end{equation}%
\begin{equation}
E_{0k}p_{k,\varepsilon }(1)=(\sigma I+A_{k})^{-2}\left( \frac{1}{\varepsilon 
}\overline{v_{k,\varepsilon }^{\ast }}(1)\right) .  \label{57-1}
\end{equation}%
Also, we introduce 
\begin{equation*}
B_{k}=F_{0k}^{\ast }E_{0k}^{-1},\mbox{ }B_{k}:D(B_{k})\subset H\rightarrow H,%
\mbox{ }D(B_{k})=\{v\in H^{2}(0,L);\mbox{ }E_{0k}^{-1}v\in
H_{0}^{2}(0,L)\}=D(A_{k}).
\end{equation*}%
By the same argument as for $A_{k}$ we infer that $-B_{k}$ generates a $%
C_{0} $-analytic semigroup.

We write the equivalent equation for $\widetilde{p_{k,\varepsilon }}%
(t):=E_{0k}p_{k,\varepsilon }(t),$ 
\begin{equation*}
-\frac{d\widetilde{p_{k,\varepsilon }}}{dt}(t)+T_{k,\varepsilon }^{\ast
}B_{k}\widetilde{p_{k,\varepsilon }}(t)=0,\mbox{ a.e. }t\in (0,1),
\end{equation*}%
\begin{equation*}
\widetilde{p_{k,\varepsilon }}(1)=(\sigma I+A_{k})^{-2}\left( \frac{1}{%
\varepsilon }\overline{v_{k,\varepsilon }^{\ast }}(1)\right) .
\end{equation*}%
The solution is 
\begin{equation*}
\widetilde{p_{k,\varepsilon }}(t)=e^{T_{k,\varepsilon }^{\ast }B_{k}(1-t)}%
\widetilde{p_{k,\varepsilon }}(1),\mbox{ for all }t\in \lbrack 0,1]
\end{equation*}%
and since $\overline{v_{k,\varepsilon }^{\ast }}(1)\in H,$ by (\ref{vk}), we
have $\widetilde{p_{k,\varepsilon }}(1)\in D(A_{k}^{2})=D(B_{k}^{2}).$
Therefore, 
\begin{equation}
\widetilde{p_{k,\varepsilon }}\in C([0,1];D(B_{k}^{2}))\cap
C^{1}([0,1];D(B_{k}))  \label{57-2-0}
\end{equation}%
and 
\begin{equation}
p_{k,\varepsilon }(t)=E_{0k}^{-1}(\widetilde{p_{k,\varepsilon }}(t)),\mbox{
for all }t\in \lbrack 0,1]  \label{57-3}
\end{equation}%
turns out to be in the spaces (\ref{57}), by the observation (\ref{iff})
made before Theorem 2.2. The proof is ended. \hfill $\square $

\medskip

\noindent Now, let us keep $w_{k,\varepsilon }^{\ast }$ fixed and give
variations to $T_{k,\varepsilon }^{\ast }.$ For $\lambda >0,$ we define $%
Z^{\lambda }=\frac{v_{k,\varepsilon }^{T_{k,\varepsilon }^{\ast }+\lambda
,w_{k,\varepsilon }^{\ast }}-v_{k,\varepsilon }^{T_{k,\varepsilon }^{\ast
},w_{k,\varepsilon }^{\ast }}}{\lambda }$ and by a straightforward
calculation we have $Z^{\lambda }\rightarrow Z$ as $\lambda \rightarrow 0,$
where $Z$ satisfies the system in variations%
\begin{eqnarray}
EZ_{t}(t)+T_{k,\varepsilon }^{\ast }FZ(t) &=&-T_{k,\varepsilon }^{\ast
}Fv_{k,\varepsilon }^{\ast }(t),\mbox{ a.e. }t\in (0,1),  \label{Z} \\
Z(t,0) &=&Z(t,L)=Z^{\prime }(t,0)=Z^{\prime }(t,L)=0,  \notag \\
Z(0,y) &=&0.  \notag
\end{eqnarray}%
Since the right hand side in the first equation in (\ref{Z}) is in $%
L^{2}(0,T_{k,\varepsilon }^{\ast };H)$ it follows that (\ref{Z}) has the
unique solution%
\begin{equation}
Z\in C([0,1];H^{2}(0,L))\cap W^{1,2}(0,1;H^{2}(0,L))\cap
L^{2}(0,1;H^{4}(0,L)).  \label{Z1}
\end{equation}

\subsection{Necessary conditions of optimality for $(P_{k,\protect%
\varepsilon })$}

We recall that $V_{1}=\{f\in H^{1}(0,1);$ $f(0)=0\}.$ Let us introduce the
set 
\begin{equation}
K_{T}=\{w\in V_{1};\mbox{ }\left\Vert w\right\Vert _{V_{1}}\leq \rho _{k}%
\sqrt{T}\}  \label{K}
\end{equation}%
and denote the normal cone to $K_{T}$ at $w$ by%
\begin{equation}
N_{K_{T}}(w)=\left\{ \chi \in V_{1}^{\ast };\mbox{ }{\rm Re}\left\langle
\chi (t),(w-w_{1})(t)\right\rangle _{V_{1}^{\ast },V_{1}}\geq 0,\mbox{ for
all }w_{1}\in K_{T}\right\} ,\mbox{ }  \label{59}
\end{equation}%
where $V_{1}^{\ast }$ is the dual of $V_{1}.$

We make a parenthesis for a discussion about this cone. It is known that 
\begin{equation*}
N_{K_{T}}(w)=\left\{ 
\begin{array}{l}
\underset{\alpha >0}{\cup }\alpha \Lambda (w),\mbox{ if }\left\Vert
w\right\Vert _{V_{1}}=\rho _{k}\sqrt{T} \\ 
\{0\},\mbox{ \ \ \ \ \ \ \ \ if }\left\Vert w\right\Vert _{V_{1}}<\rho _{k}%
\sqrt{T} \\ 
\varnothing ,\mbox{ \ \ \ \ \ \ \ \ \ \ if }\left\Vert w\right\Vert
_{V_{1}}>\rho _{k}\sqrt{T},%
\end{array}%
\right.
\end{equation*}%
where $\alpha \in \mathbb{R}_{+},$ and $\Lambda :V_{1}\rightarrow
V_{1}^{\ast }$ is the canonical isomorphism from $V_{1}$ to $V_{1}^{\ast }$.
Here, $\Lambda w=-\ddot{w}$ (see e.g., \cite{vb-springer-2010}, pp. 2-4).

By abuse of notation, we denote still by $N_{K_{T}}(w)$ the restriction of $%
N_{K_{T}}(w)$ on $L^{2}(0,1).$ In this case 
\begin{equation*}
\Lambda w=-\ddot{w},\mbox{ }\Lambda :D(\Lambda ):=\{w\in H^{2}(0,1);\mbox{ }%
w(0)=\dot{w}(1)=0\}\subset L^{2}(0,1)\rightarrow L^{2}(0,1).
\end{equation*}

Let $(T_{k,\varepsilon }^{\ast },w_{k,\varepsilon }^{\ast })$ be an optimal
controller in $(P_{k,\varepsilon }),$ with $w_{k,\varepsilon }^{\ast }\in
K_{T_{k,\varepsilon }^{\ast }}.$ For $\lambda >0,$ we set 
\begin{equation*}
w_{k,\varepsilon }^{\lambda }=w_{k,\varepsilon }^{\ast }+\lambda \omega ,%
\mbox{ where }\omega =\widetilde{w}-w_{k,\varepsilon }^{\ast },\mbox{ }%
\widetilde{w}\in K_{T_{k,\varepsilon }^{\ast }},
\end{equation*}%
that is, $\widetilde{w}\in V_{1},$ $\left\Vert \widetilde{w}\right\Vert
_{V_{1}}\leq \rho _{k}\sqrt{T_{k,\varepsilon }^{\ast }}.$ We recall the
notation 
\begin{equation}
h_{k,\varepsilon }^{\ast }(t)=\int_{0}^{t}(w_{k,\varepsilon }^{\ast
}-w_{k}^{\ast })(\tau )d\tau .  \label{51-0}
\end{equation}

\medskip

\noindent \textbf{Proposition 3.4. }\textit{Let }$(T_{k,\varepsilon }^{\ast
},w_{k,\varepsilon }^{\ast })$\textit{\ be an optimal control in }$%
(P_{k,\varepsilon })$ \textit{with the optimal state} $v_{k,\varepsilon
}^{\ast }$\textit{. Then, }%
\begin{equation}
\alpha _{k,\varepsilon }\ddot{w}_{k,\varepsilon }^{\ast
}(t)=T_{k,\varepsilon }^{\ast }\nu \overline{p_{k,\varepsilon }^{\prime
\prime \prime }}(t,L)+\int_{t}^{1}h_{k,\varepsilon }^{\ast }(\tau )d\tau ,%
\mbox{ \textit{for all} }t\in \lbrack 0,1],\mbox{ }  \label{eta}
\end{equation}%
\textit{where }$w_{k,\varepsilon }^{\ast }(0)=\dot{w}_{k,\varepsilon }(1)=0,$
$\ddot{w}_{k,\varepsilon }^{\ast }\in H^{2}(0,1),$ $\left\Vert
w_{k,\varepsilon }^{\ast }\right\Vert _{V_{1}}=\rho _{k}\sqrt{%
T_{k,\varepsilon }^{\ast }},$ $\alpha _{k,\varepsilon }\in \mathbb{R}_{+}$ 
\textit{and } 
\begin{eqnarray}
&&\alpha _{k,\varepsilon }\rho _{k}\sqrt{T_{k,\varepsilon }^{\ast }}%
\left\Vert \ddot{w}_{k,\varepsilon }^{\ast }\right\Vert _{V_{1}^{\ast
}}+T_{k,\varepsilon }^{\ast }{\rm Re}\int_{0}^{1}(\overline{%
v_{k,\varepsilon }^{\ast }}(t),F_{0k}^{\ast }p_{k,\varepsilon }(t))_{H}dt
\label{62} \\
&=&1-{\rm Re}\int_{0}^{1}\left( \int_{t}^{1}h_{\varepsilon }^{\ast }(\tau
)d\tau \right) \overline{w_{k,\varepsilon }^{\ast }}(t)dt,\mbox{ }  \notag
\end{eqnarray}%
\textit{with }$p_{k,\varepsilon }$\ \textit{the solution to the dual
backward equation }(\ref{56-1})-(\ref{56-3}).

\medskip

\noindent \textbf{Proof. }The proof will be done in two steps.

\noindent \textbf{Step 1.} \textbf{The} \textbf{first condition of optimality%
}. Let $(T_{k,\varepsilon }^{\ast },w_{k,\varepsilon }^{\ast })$ be an
optimal controller in $(P_{k,\varepsilon }).$ As already said, we keep $%
T_{k,\varepsilon }^{\ast }$ fixed and give variations to $w_{k,\varepsilon
}^{\ast }$. Using the fact that $w_{k,\varepsilon }^{\ast }$ is optimal we
can write 
\begin{equation*}
J_{k,\varepsilon }(T_{k,\varepsilon }^{\ast },w_{k,\varepsilon }^{\ast
})\leq J_{k,\varepsilon }(T,w),\mbox{ for all }w\in K_{T}.
\end{equation*}%
This holds true if we replace $T$ by $T_{k,\varepsilon }^{\ast }$ and $w$ by 
$w_{k,\varepsilon }^{\lambda }.$ By calculating 
\begin{eqnarray*}
&&\lim_{\lambda \rightarrow 0}\frac{J_{k,\varepsilon }(T_{k,\varepsilon
}^{\ast },w_{k,\varepsilon }^{\lambda })-J_{k,\varepsilon }(T_{k,\varepsilon
}^{\ast },w_{k,\varepsilon }^{\ast })}{\lambda } \\
&=&\lim_{\lambda \rightarrow 0}\frac{1}{\lambda }\left( \frac{1}{%
2\varepsilon }\left( \left\Vert (\sigma I+A_{k})^{-1}v_{k,\varepsilon
}^{T_{k,\varepsilon }^{\ast },w_{k,\varepsilon }^{\lambda }}(1)\right\Vert
_{H}^{2}-\left\Vert (\sigma I+A_{k})^{-1}v_{k,\varepsilon
}^{T_{k,\varepsilon }^{\ast },w_{k,\varepsilon }^{\ast }}(1)\right\Vert
_{H}^{2}\right) \right. \\
&&\left. +\frac{1}{2}\left( \int_{0}^{1}\left\vert
\int_{0}^{t}(w_{k,\varepsilon }^{\lambda }-w_{k}^{\ast })(\tau )d\tau
\right\vert ^{2}dt-\int_{0}^{1}\left\vert \int_{0}^{t}(w_{k,\varepsilon
}^{\ast }-w_{k}^{\ast })(\tau )d\tau \right\vert ^{2}dt\right) \right)
\end{eqnarray*}%
we obtain 
\begin{eqnarray}
&&\lim_{\lambda \rightarrow 0}\frac{J_{k,\varepsilon }(T_{k,\varepsilon
}^{\ast },w_{k,\varepsilon }^{\lambda })-J_{k,\varepsilon }(T_{k,\varepsilon
}^{\ast },w_{k,\varepsilon }^{\ast })}{\lambda }  \label{65-0} \\
&=&{\rm Re}\left\{ \int_{0}^{L}\frac{1}{\varepsilon }(\sigma I+A_{k})^{-2}%
\overline{v_{k,\varepsilon }^{\ast }}(1)\cdot
Y(1,y)dy+\int_{0}^{1}h_{k,\varepsilon }^{\ast }(t)\left( \int_{0}^{t}%
\overline{\omega }(s)ds\right) dt\right\} \geq 0.  \notag
\end{eqnarray}%
Here we used the definition of the scalar product in $\mathbb{C}$ and that $%
{\rm Re}(a\cdot \overline{b})={\rm Re}(\overline{a}\cdot b),$ for $a,b\in 
\mathbb{C}$. Then, we calculate the last term in (\ref{65-0})%
\begin{equation}
{\rm Re}\int_{0}^{1}h_{k,\varepsilon }^{\ast }(t)\left( \int_{0}^{t}%
\overline{\omega }(s)ds\right) dt={\rm Re}\int_{0}^{1}\int_{s}^{1}h_{k,%
\varepsilon }^{\ast }(t)\overline{\omega }(s)dtds  \label{g51}
\end{equation}%
and replacing it in (\ref{65-0}) we obtain%
\begin{equation}
{\rm Re}\left\{ \int_{0}^{L}(\sigma I+A_{k})^{-2}\left( \frac{1}{%
\varepsilon }\overline{v_{k,\varepsilon }^{\ast }}(1)\right)
Y(1,y)dy+\int_{0}^{1}\left( \int_{t}^{1}h_{k,\varepsilon }^{\ast }(\tau
)d\tau \right) \overline{\omega }(t)dt\right\} \geq 0.  \label{65}
\end{equation}%
Using the notation (\ref{L-V}) we can write eq. (\ref{54}) in the equivalent
form%
\begin{equation}
EY_{t}(t)+T_{k,\varepsilon }^{\ast }FY(t)=0  \label{65-1}
\end{equation}%
with the boundary and initial conditions. Now, we multiply this equation
scalarly in $H$ by $p_{k,\varepsilon }(t)$ and integrate with respect to $t$
over $(0,1).$ While performing all the integrals by parts, we put first into
evidence the relation 
\begin{equation}
\int_{0}^{L}FY(t,y)\cdot p_{k,\varepsilon }(t,y)dy=\int_{0}^{L}F_{0k}^{\ast
}p_{k,\varepsilon }(t,y)\cdot Y(t,y)dy-\nu p_{k,\varepsilon }^{\prime \prime
\prime }(t,L)\omega (t),  \label{F-F0}
\end{equation}%
and obtain%
\begin{eqnarray}
&&\int_{0}^{1}\int_{0}^{L}\left\{ -(k^{2}p_{k,\varepsilon }-p_{k,\varepsilon
}^{\prime \prime })_{t}\right\} Ydydt  \label{66} \\
&&+T_{k,\varepsilon }^{\ast }\int_{0}^{1}\int_{0}^{L}\{\nu p_{k,\varepsilon
}^{\mbox{iv}}-(2\nu k^{2}+ikU)p_{k,\varepsilon }^{\prime \prime
}-2ikU^{\prime }p_{k,\varepsilon }^{\prime }+(\nu
k^{4}+ik^{3}U)p_{k,\varepsilon }\}Ydydt  \notag \\
&&+\int_{0}^{L}\left\{ k^{2}p_{k,\varepsilon }(1,y)-p_{k,\varepsilon
}^{\prime \prime }(1,y)\right\} Y(1,y)dy-T_{k,\varepsilon }^{\ast
}\int_{0}^{1}\omega (t)\nu p_{k,\varepsilon }^{\prime \prime \prime
}(t,L)dy=0.  \notag
\end{eqnarray}%
By (\ref{57}), we see that $p_{k,\varepsilon }^{\prime \prime \prime }\in
C([0,1];H^{1}(0,L))$ and so\ the trace $p_{k,\varepsilon }^{\prime \prime
\prime }(t,L)$ is well defined and belongs to $C([0,1];\mathbb{C}).$
Recalling now the equations in the adjoint system (\ref{56-1})-(\ref{56-3})
and comparing with the integrands in (\ref{66}) we get%
\begin{equation}
\int_{0}^{L}(\sigma I+A_{k})^{-2}\left( \frac{1}{\varepsilon }\overline{%
v_{k,\varepsilon }^{\ast }}(1)\right) \cdot Y(1,y)dy=T_{k,\varepsilon
}^{\ast }\int_{0}^{1}\omega (t)\cdot \nu p_{k,\varepsilon }^{\prime \prime
\prime }(t,L)dt.  \label{67}
\end{equation}%
Then plugging the latter into (\ref{65}) and using again ${\rm Re}(a\cdot
b)={\rm Re}(\overline{a}\cdot \overline{b})$ we deduce the relation 
\begin{eqnarray}
&&{\rm Re}\left\{ \int_{0}^{1}T_{k,\varepsilon }^{\ast }\nu \overline{%
p_{k,\varepsilon }^{\prime \prime \prime }}(t,L)\overline{\omega }%
(t)dt+\int_{0}^{1}\left( \int_{t}^{1}h_{k,\varepsilon }^{\ast }(\tau )d\tau
\right) \overline{\omega }(t)dt\right\}  \label{68} \\
&=&{\rm Re}\int_{0}^{1}\left( -T_{k,\varepsilon }^{\ast }\nu \overline{%
p_{\varepsilon }^{\prime \prime \prime }}(t,L)-\int_{t}^{1}h_{\varepsilon
}^{\ast }(\tau )d\tau \right) (w_{k,\varepsilon }^{\ast }-\widetilde{w}%
)(t)dt\geq 0,  \notag
\end{eqnarray}%
for all $\widetilde{w}\in K_{T_{k,\varepsilon }^{\ast }}$. The immediate
result is that 
\begin{equation}
\eta _{k,\varepsilon }:=-T_{k,\varepsilon }^{\ast }\nu \overline{%
p_{\varepsilon }^{\prime \prime \prime }}(\cdot ,L)-\int_{\cdot
}^{1}h_{\varepsilon }^{\ast }(\tau )d\tau \in N_{K_{T_{k,\varepsilon }^{\ast
}}}(w_{k,\varepsilon }^{\ast })  \label{68-0}
\end{equation}%
where $N_{K_{T_{k,\varepsilon }^{\ast }}}(w_{k,\varepsilon }^{\ast })$ is
viewed as the cone from $V_{1}$ to $V_{1}^{\ast }.$ However, since $\eta
_{k,\varepsilon }\in L^{2}(0,1),$ we can consider it in the restriction of $%
N_{K_{T_{k,\varepsilon }^{\ast }}}(w_{k,\varepsilon }^{\ast })$ to $%
L^{2}(0,1),$ still denoted by $N_{K_{T_{k,\varepsilon }^{\ast
}}}(w_{k,\varepsilon }^{\ast }).$ Because an element of this cone is of the
form $\eta _{k,\varepsilon }=-\alpha _{k,\varepsilon }\ddot{w}%
_{k,\varepsilon }$ with $w_{k,\varepsilon }^{\ast }(0)=\dot{w}%
_{k,\varepsilon }^{\ast }(1)=0$ and $\alpha _{k,\varepsilon }>0$, we have 
\begin{equation}
-\alpha _{k,\varepsilon }\ddot{w}_{k,\varepsilon }(t)=-T_{k,\varepsilon
}^{\ast }\nu \overline{p_{k,\varepsilon }^{\prime \prime \prime }}%
(t,L)-\int_{t}^{1}h_{k,\varepsilon }^{\ast }(\tau )d\tau ,\mbox{ for all }%
t\in \lbrack 0,1],  \label{70-0}
\end{equation}%
which is just (\ref{eta}). We note that $\ddot{w}_{k,\varepsilon }$ turns
out to be in $H^{2}(0,1).$

The situation in which $\left\Vert w_{k,\varepsilon }^{\ast }\right\Vert
_{V_{1}}<\rho _{k}\sqrt{T_{k,\varepsilon }^{\ast }}$ provides $\ddot{w}%
_{k,\varepsilon }(t)=0,$ which gives the solution $w_{k,\varepsilon }^{\ast
}(t)=0$ which is not relevant for our problem. In particular, it means that $%
w_{k}^{\ast }=0,$ that is the flow would be not controlled.

Thus, (\ref{eta}) follows for the case when $\left\Vert w_{k,\varepsilon
}^{\ast }\right\Vert _{V_{1}}=\rho _{k}\sqrt{T_{k,\varepsilon }^{\ast }}.$

\textbf{Step 2. The second condition of optimality. }Here\textbf{\ }we keep $%
w_{k,\varepsilon }^{\ast }$ fixed and give variations for $T_{k,\varepsilon
}^{\ast }.$ Since $T_{k,\varepsilon }^{\ast }$ realizes the minimum in $%
(P_{k,\varepsilon })$ we write 
\begin{equation*}
J_{k,\varepsilon }(T_{k,\varepsilon }^{\ast },w_{k,\varepsilon }^{\ast
})\leq J_{k,\varepsilon }(T_{k,\varepsilon }^{\ast }+\lambda
,w_{k,\varepsilon }^{\ast }),\mbox{ }\lambda >0,
\end{equation*}%
that is, 
\begin{eqnarray*}
&&J_{k,\varepsilon }(T_{k,\varepsilon }^{\ast },w_{k,\varepsilon }^{\ast
})=T_{k,\varepsilon }^{\ast }+\frac{1}{2\varepsilon }\left\Vert (\sigma
I+A_{k})^{-2}v_{k,\varepsilon }^{T_{k,\varepsilon }^{\ast },w_{k,\varepsilon
}^{\ast }}(1)\right\Vert _{H}^{2}+\frac{1}{2}\int_{0}^{1}\left\vert
h_{k,\varepsilon }^{\ast }(t)\right\vert ^{2}dt \\
&\leq &J_{k,\varepsilon }(T_{k,\varepsilon }^{\ast }+\lambda
,w_{k,\varepsilon }^{\ast })=T_{k,\varepsilon }^{\ast }+\lambda +\frac{1}{%
2\varepsilon }\left\Vert (\sigma I+A_{k})^{-2}v_{k,\varepsilon
}^{T_{k,\varepsilon }^{\ast }+\lambda ,w_{k,\varepsilon }^{\ast
}}(1)\right\Vert _{H}^{2}+\frac{1}{2}\int_{0}^{1}\left\vert h_{k,\varepsilon
}^{\ast }(t)\right\vert ^{2}dt.
\end{eqnarray*}%
By performing the computations as before we obtain 
\begin{equation}
1+{\rm Re}\int_{0}^{1}\frac{1}{\varepsilon }(\sigma
I+A_{k})^{-2}v_{k,\varepsilon }^{\ast }(1)Z(1,y)dy=0,  \label{cond2}
\end{equation}%
where $Z$ is the solution to the system in variations (\ref{Z}). We consider
the same adjoint system (\ref{57-0})-(\ref{57-1}). By multiplying the first
equation in (\ref{Z}) by $p_{k,\varepsilon }(t)$ scalarly in $H,$
integrating along $t\in (0,1)$ and taking the real part, we obtain 
\begin{equation}
{\rm Re}\int_{0}^{L}E_{0k}p_{k,\varepsilon }(1)Z(1,y)dy=-{\rm Re}%
\int_{0}^{1}\int_{0}^{L}T_{k,\varepsilon }^{\ast }Fv_{k,\varepsilon }^{\ast
}p_{k,\varepsilon }dydt.  \label{cond2-1}
\end{equation}%
By using the final condition in the adjoint system and comparing (\ref{cond2}%
) and (\ref{cond2-1}) we get%
\begin{equation}
{\rm Re}\int_{0}^{1}\int_{0}^{L}T_{k,\varepsilon }^{\ast }Fv_{k,\varepsilon
}^{\ast }p_{k,\varepsilon }dydt=1.  \label{cond2-3}
\end{equation}%
Finally, applying a relation similar to (\ref{F-F0}) in the left-hand side
of (\ref{cond2-3}) we obtain%
\begin{equation*}
{\rm Re}\int_{0}^{1}\int_{0}^{L}T_{k,\varepsilon }^{\ast }F_{0k}^{\ast
}p_{k,\varepsilon }\cdot v_{k,\varepsilon }^{\ast }dydt-{\rm Re}%
\int_{0}^{1}T_{k,\varepsilon }^{\ast }\nu p_{k,\varepsilon }^{\prime \prime
\prime }(t,L)w_{k,\varepsilon }^{\ast }dt=1
\end{equation*}%
and recalling again that ${\rm Re}(a\cdot b)={\rm Re}(\overline{a}\cdot 
\overline{b})$ we get%
\begin{equation}
-{\rm Re}\int_{0}^{1}T_{k,\varepsilon }^{\ast }\nu \overline{%
p_{k,\varepsilon }^{\prime \prime \prime }}(t,L)\overline{w_{k,\varepsilon
}^{\ast }}dt+{\rm Re}\int_{0}^{1}T_{k,\varepsilon }^{\ast }(\overline{%
v_{k,\varepsilon }^{\ast }}(t),F_{0k}^{\ast }p_{k,\varepsilon }(t))_{H}dt=1.
\label{cond2-2}
\end{equation}%
Using (\ref{eta}) we still can write 
\begin{eqnarray}
&&{\rm Re}\int_{0}^{1}\left( -T_{k,\varepsilon }^{\ast }\nu
p_{k,\varepsilon }^{\prime \prime \prime }(t,L)-\int_{t}^{1}h_{\varepsilon
}^{\ast }(\tau )d\tau \right) \overline{w_{k,\varepsilon }^{\ast }}(t)dt+%
{\rm Re}\int_{0}^{1}T_{k,\varepsilon }^{\ast }(\overline{v_{k,\varepsilon
}^{\ast }}(t),F_{0k}^{\ast }p_{k,\varepsilon }(t))_{H}dt  \notag \\
&=&1-{\rm Re}\int_{0}^{1}\left( \int_{t}^{1}h_{\varepsilon }^{\ast }(\tau
)d\tau \right) \overline{w_{k,\varepsilon }^{\ast }}dt.  \label{cond3}
\end{eqnarray}%
But $\eta _{k,\varepsilon }(t)\in N_{K_{T_{k,\varepsilon }^{\ast
}}}(w_{k,\varepsilon }^{\ast })$ and so we have 
\begin{equation}
{\rm Re}\int_{0}^{1}\eta _{k,\varepsilon }(t)\overline{w_{k,\varepsilon
}^{\ast }}(t)dt={\rm Re}\left\langle \eta _{k,\varepsilon
}(t),w_{k,\varepsilon }^{\ast }(t)\right\rangle _{V_{1}^{\ast },V_{1}}=\rho
_{k}\sqrt{T_{k,\varepsilon }^{\ast }}\left\Vert \eta _{k,\varepsilon
}\right\Vert _{V_{1}^{\ast }}.  \label{cond4}
\end{equation}%
Let $\left\Vert w_{k,\varepsilon }^{\ast }\right\Vert _{V_{1}}=\rho _{k}%
\sqrt{T_{k,\varepsilon }^{\ast }}.$ Taking into account that in this case
there is $\alpha _{k,\varepsilon }>0$ such that $\eta _{k,\varepsilon
}(t)=-\alpha _{k,\varepsilon }\ddot{w}_{k,\varepsilon }^{\ast },$ we get by (%
\ref{cond3}) and (\ref{cond4}), relation (\ref{62}), as claimed.

Moreover, on the one hand, $\left\langle -\alpha _{k,\varepsilon }\ddot{w}%
_{k,\varepsilon }^{\ast },\overline{w_{k,\varepsilon }^{\ast }}%
(t)\right\rangle _{V_{1}^{\ast },V_{1}}=\alpha _{k,\varepsilon
}\int_{0}^{1}\left\vert \dot{w}_{k,\varepsilon }^{\ast }(t)\right\vert
^{2}dt=\alpha _{k,\varepsilon }\rho _{k}^{2}T_{k,\varepsilon }^{\ast }.$ On
the other hand by (\ref{cond4}) we get 
\begin{equation}
{\rm Re}\left\langle -\alpha _{k,\varepsilon }\ddot{w}_{k,\varepsilon
}^{\ast },\overline{w_{k,\varepsilon }^{\ast }}(t)\right\rangle
_{V_{1}^{\ast },V_{1}}=\alpha _{k,\varepsilon }\rho _{k}\sqrt{%
T_{k,\varepsilon }^{\ast }}\left\Vert \ddot{w}_{k,\varepsilon }^{\ast
}\right\Vert _{V_{1}^{\ast }},  \label{cond5}
\end{equation}%
hence $\left\Vert \ddot{w}_{k,\varepsilon }^{\ast }\right\Vert _{V_{1}^{\ast
}}=\rho _{k}\sqrt{T_{k,\varepsilon }^{\ast }}$ is verified.

The situation $\left\Vert \ddot{w}_{k,\varepsilon }^{\ast }\right\Vert
_{V_{1}^{\ast }}<\rho _{k}\sqrt{T_{k,\varepsilon }^{\ast }}$ was excluded,
so that (\ref{eta}) and (\ref{62}) follow for the case when $\left\Vert 
\ddot{w}_{k,\varepsilon }^{\ast }\right\Vert _{V_{1}^{\ast }}=\rho _{k}\sqrt{%
T_{k,\varepsilon }^{\ast }}.$ This ends the proof. \hfill $\square $

\section{The maximum principle for $(P_{k})$}

\setcounter{equation}{0}

We recall that by $(\widehat{H_{k}}),$ which is $(H_{k})$ translated in the
new variable (\ref{t-cap}), the system (\ref{4'})-(\ref{4''''}) is
controllable. Namely, for each initial datum $v^{t_{0}}\in L^{2}(0,L),$ $%
\left\Vert v^{t_{0}}\right\Vert _{(H^{2}(0,L))^{\ast }}\leq 1,$ there exists 
$w\in H^{1}(0,1),$ and $\gamma _{(\cdot ,1)}\in L^{2}(0,1),$ satisfying $%
w(\tau )=0$ for $0\leq \tau \leq t_{0},$ $\left( \int_{0}^{1}\left\vert \dot{%
w}(\tau )\right\vert ^{2}d\tau \right) ^{1/2}\leq \gamma _{(t_{0},1)},$ such
that $v^{t_{0},w}(1,y)=0$ a.e. $y\in (0,L).$ Here, $v^{t_{0},w}$ is the
solution to (\ref{4'})-(\ref{4''''}) starting at time $t_{0}$. (

\noindent Here, all functions are in fact those denoted by $"$ $\widehat{}$ $%
"$ depending on $\widehat{t},$ but this decoration is skipped in this
section).

We denote by $C_{\nu ,U,k}:=Ck^{4}\left( \nu +\left\Vert U\right\Vert
_{W^{1,\infty }(0,L)}\right) ,$ where $C$ is a positive number$.$

\medskip

\noindent \textbf{Theorem 4.1. }\textit{Let }$(T_{k}^{\ast },w_{k}^{\ast
},v_{k}^{\ast })$\textit{\ be optimal in }$(\widehat{P_{k}})$\textit{. If }%
\begin{eqnarray}
&&C_{\nu ,U,k}(1+k^{2})\left\Vert \gamma \right\Vert _{L^{1}(0,1)}\sqrt{%
T_{k}^{\ast }}\left( T_{k}^{\ast }+\sqrt{T_{k}^{\ast }}\right) <1,\mbox{ }
\label{cond-rok} \\
&&\rho _{k}\left( 1-C_{\nu ,U,k}(1+k^{2})\left\Vert \gamma \right\Vert
_{L^{1}(0,1)}\sqrt{T_{k}^{\ast }}(T_{k}^{\ast }+\sqrt{T_{k}^{\ast }})\right)
\notag \\
&>&C_{\nu ,U,k}(1+k^{2})\sqrt{T_{k,\varepsilon }^{\ast }}\left\Vert \gamma
\right\Vert _{L^{1}(0,1)}\left\Vert v_{k0}\right\Vert _{H^{4}(0,L)\cap
H_{0}^{2}(0,L)},  \notag
\end{eqnarray}%
\textit{then, there exists }$\alpha _{k}^{\ast }>0$ \textit{such that }%
\begin{equation}
\alpha _{k}^{\ast }\ddot{w}_{k}^{\ast }(t)=T_{k}^{\ast }\nu \overline{%
p_{k}^{\prime \prime \prime }}(t,L),\mbox{ \textit{a.e.} }t\in (0,1),\mbox{ }%
w_{k}^{\ast }(0)=\dot{w}_{k}^{\ast }(1)=0,\mbox{ }\left\Vert w_{k}^{\ast
}\right\Vert _{V_{1}}=\rho _{k}\sqrt{T_{k}^{\ast }},  \label{101}
\end{equation}%
\textit{\ } 
\begin{equation}
\alpha _{k}^{\ast }\rho _{k}^{2}T_{k}^{\ast }+T_{k}^{\ast }{\rm Re}%
\int_{0}^{1}(\overline{v_{k}^{\ast }}(t),F_{0k}^{\ast }p_{k}(t))_{H}dt=1,
\label{102}
\end{equation}%
\textit{where }$p_{k}$ \textit{is the solution to the adjoint equation}%
\begin{equation}
-E_{0k}(p_{k})_{t}(t)+T_{k}^{\ast }F_{0k}^{\ast }p_{k}(t)=0,\mbox{ \textit{%
a.e.} }t\in (0,1).  \label{103}
\end{equation}

\medskip

\noindent \textbf{Proof. }Using the notation (\ref{L-V}) let us consider the
generic system 
\begin{eqnarray}
Ev_{t}(t)+T_{k,\varepsilon }^{\ast }Fv(t) &=&0,  \notag \\
v(t,0) &=&0,\mbox{ }v(t,L)=w(t),\mbox{ }v^{\prime }(t,0)=v^{\prime }(t,L)=0,
\label{91} \\
v(0) &=&v_{0k}\in L^{2}(0,L),  \notag
\end{eqnarray}%
which has a unique solution, according to Theorem 2.2 and which is
controllable, according to $(\widehat{H_{k}}).$ Then, we repeat a similar
calculus providing (\ref{66}). Thus, we multiply the first equation in (\ref%
{91}) scalarly by $p_{k,\varepsilon }(t),$ integrate along $(t,1)$ and get%
\begin{eqnarray}
&&\int_{0}^{L}v(1,y)E_{0k}p_{k,\varepsilon
}(1)dy-\int_{0}^{L}v(t,y)E_{0k}p_{k,\varepsilon }(t)dy  \label{91-0} \\
&&-\int_{t}^{1}T_{k,\varepsilon }^{\ast }\nu p_{k,\varepsilon }^{\prime
\prime \prime }(\tau ,L)w(\tau )d\tau =0,\mbox{ for }t\in (0,1).  \notag
\end{eqnarray}%
Now, we apply $(\widehat{H_{k}})$ for a particular choice, for $t_{0}=t,$ $%
v^{t_{0}}=v^{t}:=v(t),$ $\left\Vert v(t)\right\Vert _{(H^{2}(0,L))^{\ast
}}\leq 1,$ $w\in V_{1},$ $w(\tau )=0$ for $0\leq \tau \leq t,$ $\gamma
_{(t,1)}:=\gamma (t)$ with $\gamma \in L^{2}(0,1)$ and $\left\Vert
w\right\Vert _{V_{1}}\leq \gamma (t).$ Then, the solution $v^{t,w}$ to (\ref%
{91}), corresponding to the initial datum $v^{t,w}(t,y)=v(t,y),$ satisfies $%
v^{t,w}(1,y)=0.$ By taking the real part in (\ref{91-0}) we get%
\begin{equation}
{\rm Re}(E_{0k}p_{k,\varepsilon }(t),\overline{v}(t))_{H}=-\nu
T_{k,\varepsilon }^{\ast }{\rm Re}\int_{t}^{1}\overline{p_{k,\varepsilon
}^{\prime \prime \prime }}(\tau ,L)\overline{w}(\tau )d\tau .  \label{91-00}
\end{equation}

We choose now $\overline{v}(t)=\frac{E_{0k}p_{k,\varepsilon }(t)}{\left\Vert
E_{0k}p_{k,\varepsilon }(t)\right\Vert _{H}}$ which ensures that $\left\Vert 
\overline{v}(t)\right\Vert _{H}=1,$ so that $\left\Vert \overline{v}%
(t)\right\Vert _{(H^{2}(0,L))^{\ast }}\leq 1$ and compute 
\begin{equation*}
\left\Vert E_{0k}p_{k,\varepsilon }(t)\right\Vert _{H}^{2}=k^{4}\left\Vert
p_{k,\varepsilon }^{\prime \prime }(t)\right\Vert _{H}^{2}+2k^{2}\left\Vert
p_{k,\varepsilon }^{\prime }(t)\right\Vert _{H}^{2}+\left\Vert
p_{k,\varepsilon }(t)\right\Vert _{H}^{2}\geq 2\left\Vert p_{k,\varepsilon
}(t)\right\Vert _{H^{2}(0,L)}^{2},
\end{equation*}%
for $k\geq 1$.

Since $\overline{w}(\tau )=0$ for $\tau \in \lbrack 0,t],$ $p_{k,\varepsilon
}^{\prime \prime \prime }(\cdot ,L)\in L^{2}(0,1)$ and $V_{1}\subset
H^{1}(0,1)\subset L^{2}(0,1)\subset (H^{1}(0,1))^{\ast }\subset V_{1}^{\ast
},$ we have by (\ref{91-00})%
\begin{eqnarray}
\left\Vert p_{k,\varepsilon }(t)\right\Vert _{H^{2}(0,L)} &\leq &\frac{1}{%
\sqrt{2}}\left\Vert E_{0k}p_{k,\varepsilon }(t)\right\Vert _{H}\leq \frac{1}{%
\sqrt{2}}\nu T_{k,\varepsilon }^{\ast }\left\vert \int_{0}^{1}\overline{%
p_{k,\varepsilon }^{\prime \prime \prime }}(\tau ,L)\overline{w}(\tau )d\tau
\right\vert  \label{201} \\
&\leq &\frac{1}{\sqrt{2}}\nu T_{k,\varepsilon }^{\ast }\left\vert
\left\langle \overline{p_{k,\varepsilon }^{\prime \prime \prime }}(\tau ,L),%
\overline{w}(\tau )\right\rangle _{V_{1}^{\ast },V_{1}}\right\vert \leq 
\frac{1}{\sqrt{2}}\nu T_{k,\varepsilon }^{\ast }\left\Vert \overline{%
p_{k,\varepsilon }^{\prime \prime \prime }}(\cdot ,L)\right\Vert
_{V_{1}^{\ast }}\left\Vert \overline{w}\right\Vert _{V_{1}}  \notag \\
&\leq &\frac{1}{\sqrt{2}}\nu T_{k,\varepsilon }^{\ast }\left\Vert \overline{%
p_{k,\varepsilon }^{\prime \prime \prime }}(\cdot ,L)\right\Vert
_{V_{1}^{\ast }}\gamma (t),  \notag
\end{eqnarray}%
where we took into account that $\left\Vert w\right\Vert _{V_{1}}\leq \gamma
(t),$ by $(\widehat{H_{k}}).$ Net, we note that 
\begin{equation}
\left\Vert \psi \right\Vert _{H^{2}(0,L)}=\sup_{\left\Vert \omega
\right\Vert _{(H^{2}(0,L))^{\ast }}\leq 1}\left\vert \left\langle \psi
,\omega \right\rangle _{H^{2}(0,L),(H^{2}(0,L))^{\ast }}\right\vert .
\label{201-0}
\end{equation}%
Going back to (\ref{91-00}) 
\begin{equation*}
{\rm Re}\left\langle E_{0k}p_{k,\varepsilon }(t),\overline{v}%
(t)\right\rangle _{H^{2}(0,L),(H^{2}(0,L))^{\ast }}={\rm Re}%
(E_{0k}p_{k,\varepsilon }(t),\overline{v}(t))_{H}=-\nu T_{k,\varepsilon
}^{\ast }{\rm Re}\int_{t}^{1}\overline{p_{k,\varepsilon }^{\prime \prime
\prime }}(\tau ,L)\overline{w}(\tau )d\tau
\end{equation*}%
and taking $\sup $, we deduce by (\ref{201-0}) that 
\begin{eqnarray}
\left\Vert E_{0k}p_{k,\varepsilon }(t)\right\Vert _{H^{2}(0,L)}
&=&\sup_{\left\Vert v(t)\right\Vert _{(H^{2}(0,L)^{\ast }}\leq 1}\left\vert
\left\langle E_{0k}p_{k,\varepsilon }(t),\overline{v}(t)\right\rangle
_{H^{2}(0,L),(H^{2}(0,L))^{\ast }}\right\vert  \label{200} \\
&\leq &\nu T_{k,\varepsilon }^{\ast }\left\Vert \overline{p_{k,\varepsilon
}^{\prime \prime \prime }}(\cdot ,L)\right\Vert _{V_{1}^{\ast }}\left\Vert 
\overline{w}\right\Vert _{V_{1}}\leq \nu T_{k,\varepsilon }^{\ast
}\left\Vert \overline{p_{k,\varepsilon }^{\prime \prime \prime }}(\cdot
,L)\right\Vert _{V_{1}^{\ast }}\gamma (t),  \notag
\end{eqnarray}%
where we used again the part of $(\widehat{H_{k}})$ written before. Now, by
a straightforward calculation we deduce that 
\begin{equation}
\left\Vert p_{k,\varepsilon }(t)\right\Vert _{H^{4}(0,L)}\leq \sqrt{2}%
\left\Vert E_{0k}p_{k,\varepsilon }(t)\right\Vert _{H^{2}(0,L)}+\sqrt{2}%
k^{2}\left\Vert p_{k,\varepsilon }^{\prime \prime }(t)\right\Vert
_{L^{2}(0,L)}.  \label{94-0}
\end{equation}%
Using now (\ref{200}) and (\ref{201}) we write 
\begin{equation*}
\left\Vert p_{k,\varepsilon }(t)\right\Vert _{H^{4}(0,L)}\leq \sqrt{2}\nu
T_{k,\varepsilon }^{\ast }\left\Vert \overline{p_{k,\varepsilon }^{\prime
\prime \prime }}(\cdot ,L)\right\Vert _{V_{1}^{\ast }}\gamma (t)+k^{2}\nu
T_{k,\varepsilon }^{\ast }\left\Vert \overline{p_{k,\varepsilon }^{\prime
\prime \prime }}(\cdot ,L)\right\Vert _{V_{1}^{\ast }}\gamma (t)
\end{equation*}%
and obtain the observability relation%
\begin{equation}
\left\Vert p_{k,\varepsilon }(t)\right\Vert _{H^{4}(0,L)}\leq (\sqrt{2}%
+k^{2})\nu T_{k,\varepsilon }^{\ast }\left\Vert \overline{p_{k,\varepsilon
}^{\prime \prime \prime }}(\cdot ,L)\right\Vert _{V_{1}^{\ast }}\gamma (t).
\label{94}
\end{equation}%
By (\ref{57-2}) we have the relation $\left\Vert F_{0k}^{\ast
}p_{k,\varepsilon }(t)\right\Vert _{H}\leq C_{\nu ,U,k}\left\Vert
p_{k,\varepsilon }(t)\right\Vert _{H^{4}(0,L)}.$

Now, we rewrite (\ref{eta}) (multiplied by $\rho _{k}\sqrt{T_{k,\varepsilon
}^{\ast }})$ and replace $\rho _{k}\sqrt{T_{k,\varepsilon }^{\ast }}\alpha
_{k,\varepsilon }\left\Vert \ddot{w}_{k,\varepsilon }^{\ast }\right\Vert
_{V_{1}^{\ast }}$ using (\ref{62}), as follows:%
\begin{eqnarray*}
&&\rho _{k}\sqrt{T_{k,\varepsilon }^{\ast }}\left\Vert T_{k,\varepsilon
}^{\ast }\nu \overline{p_{k,\varepsilon }^{\prime \prime \prime }}(\cdot
,L)\right\Vert _{V_{1}^{\ast }}\leq \rho _{k}\sqrt{T_{k,\varepsilon }^{\ast }%
}\alpha _{k,\varepsilon }\left\Vert \ddot{w}_{k,\varepsilon }^{\ast
}\right\Vert _{V_{1}^{\ast }}+\rho _{k}\sqrt{T_{k,\varepsilon }^{\ast }}%
\left\Vert \int_{\cdot }^{1}h_{k,\varepsilon }^{\ast }(\tau )d\tau
\right\Vert _{V_{1}^{\ast }} \\
&\leq &T_{k,\varepsilon }^{\ast }\int_{0}^{1}\left\Vert \overline{%
v_{k,\varepsilon }^{\ast }}(t)\right\Vert _{H}\left\Vert F_{0k}^{\ast
}p_{k,\varepsilon }(t)\right\Vert _{H}dt+1+\int_{0}^{1}\left\vert
\int_{t}^{1}h_{k,\varepsilon }^{\ast }(\tau )d\tau \right\vert \left\vert 
\overline{w_{k,\varepsilon }^{\ast }}(t)\right\vert dt+C_{\varepsilon }(\rho
_{k},T_{k,\varepsilon }^{\ast }) \\
&\leq &T_{k,\varepsilon }^{\ast }C_{v}\int_{0}^{1}C_{\nu ,U,k}\left\Vert
p_{k,\varepsilon }(t)\right\Vert _{H^{4}(0,L)}dt+1+C_{\varepsilon }(\rho
_{k},T_{k,\varepsilon }^{\ast }) \\
&\leq &C_{v}C_{\nu ,U,k}T_{k,\varepsilon }^{\ast }(\sqrt{2}%
+k^{2})\int_{0}^{1}\nu T_{k,\varepsilon }^{\ast }\left\Vert \overline{%
p_{k,\varepsilon }^{\prime \prime \prime }}(\cdot ,L)\right\Vert
_{V_{1}^{\ast }}\gamma (t)dt+1+C_{\varepsilon }(\rho _{k},T_{k}^{\ast }) \\
&\leq &C_{v}C_{\nu ,U,k}T_{k}^{\ast }(\sqrt{2}+k^{2})\nu T_{k}^{\ast
}\left\Vert \overline{p_{k,\varepsilon }^{\prime \prime \prime }}(\cdot
,L)\right\Vert _{V_{1}^{\ast }}\left\Vert \gamma \right\Vert
_{L^{1}(0,1)}+1+C_{\varepsilon }(\rho _{k},T_{k}^{\ast }),
\end{eqnarray*}%
where $C_{\varepsilon }(\rho _{k},T_{k}^{\ast })$ goes to zero, as $%
\varepsilon \rightarrow 0,$ by (\ref{52-6}).

\noindent Here, we took into account the observability relation (\ref{94}),
the boundedness of $v_{k,\varepsilon }^{\ast }$ in $C([0,1];H^{2}(0,L))$ by (%
\ref{vk-est1}) (with the constant in (\ref{vk-est1}) called here $C_{v})$
and the fact that $T_{k,\varepsilon }^{\ast }\leq T_{k}^{\ast }$ by (\ref%
{52-1}). Finally, we obtain%
\begin{eqnarray}
&&\rho _{k}\sqrt{T_{k,\varepsilon }^{\ast }}\nu T_{k,\varepsilon }^{\ast
}\left\Vert \overline{p_{k,\varepsilon }^{\prime \prime \prime }}(\cdot
,L)\right\Vert _{V_{1}^{\ast }}\leq 1+C_{\varepsilon }(\rho
_{k},T_{k,\varepsilon }^{\ast })+\sqrt{2}C_{\nu ,U,k}T_{k}^{\ast
}(1+k^{2})\left\Vert \gamma \right\Vert _{L^{1}(0,1)}  \label{95} \\
&&\times \sqrt{C}\left( \left\Vert v_{k0}\right\Vert _{H^{4}(0,L)\cap
H_{0}^{2}(0,L)}+\rho _{k}(T_{k,\varepsilon }^{\ast }+\sqrt{T_{k,\varepsilon
}^{\ast }})\right) \nu T_{k,\varepsilon }^{\ast }\left\Vert \overline{%
p_{k,\varepsilon }^{\prime \prime \prime }}(\cdot ,L)\right\Vert
_{V_{1}^{\ast }}.  \notag
\end{eqnarray}%
The constants $\sqrt{2}$ and $\sqrt{C}$ coming from (\ref{vk-est1}) will be
included in $C_{\nu ,U,k}.$ Therefore, by denoting 
\begin{eqnarray*}
D_{k}(T) &:&=1-C_{\nu ,U,k}(1+k^{2})\left\Vert \gamma \right\Vert
_{L^{1}(0,1)}\sqrt{T}(T+\sqrt{T}),\mbox{ } \\
G_{k}(T) &:&=C_{\nu ,U,k}(1+k^{2})\sqrt{T}\left\Vert \gamma \right\Vert
_{L^{1}(0,1)}\left\Vert v_{k0}\right\Vert _{H^{4}(0,L)\cap H_{0}^{2}(0,L)},
\end{eqnarray*}%
we get%
\begin{equation}
T_{k,\varepsilon }^{\ast }\nu \left\Vert \overline{p_{k,\varepsilon
}^{\prime \prime \prime }}(\cdot ,L)\right\Vert _{V_{1}^{\ast }}\sqrt{%
T_{k,\varepsilon }^{\ast }}(\rho _{k}D_{k}(T_{k,\varepsilon }^{\ast
})-G_{k}(T_{k,\varepsilon }^{\ast }))\leq 1+C_{\varepsilon }(\rho
_{k},T_{k}^{\ast }).  \label{96}
\end{equation}%
On the one hand, on the basis of (\ref{cond-rok}) we have that 
\begin{equation*}
1>C_{\nu ,U,k}(1+k^{2})\sqrt{T_{k}^{\ast }}\left( T_{k}^{\ast }+\sqrt{%
T_{k}^{\ast }}\right) \left\Vert \gamma \right\Vert _{L^{1}(0,1)}\geq C_{\nu
,U,k}(1+k^{2})\sqrt{T_{k,\varepsilon }^{\ast }}\left( T_{k,\varepsilon
}^{\ast }+\sqrt{T_{k,\varepsilon }^{\ast }}\right) \left\Vert \gamma
\right\Vert _{L^{1}(0,1)},
\end{equation*}%
since $T_{k,\varepsilon }^{\ast }<T_{k}^{\ast }$ by Theorem 3.2, so $%
D_{k}(T_{k,\varepsilon }^{\ast })>0$. On the other hand, $G_{k}(T_{k}^{\ast
})>G_{k}(T_{k,\varepsilon }^{\ast })$ and $D(T_{k}^{\ast
})<D(T_{k,\varepsilon }^{\ast }),$ so that $\rho _{k}D_{k}(T_{k,\varepsilon
}^{\ast })-G_{k}(T_{k,\varepsilon }^{\ast })>\rho _{k}D_{k}(T_{k}^{\ast
})-G_{k}(T_{k}^{\ast }).$ All these imply%
\begin{eqnarray}
\left\Vert \overline{p_{k,\varepsilon }^{\prime \prime \prime }}(\cdot
,L)\right\Vert _{V_{1}^{\ast }} &\leq &\frac{1+C_{\varepsilon }(\rho
_{k},T_{k}^{\ast })}{\nu T_{k,\varepsilon }^{\ast }\sqrt{T_{k,\varepsilon
}^{\ast }}(\rho _{k}D_{k}(T_{k,\varepsilon }^{\ast })-G_{k}(T_{k,\varepsilon
}^{\ast }))}  \label{97} \\
&<&\frac{1+C_{\varepsilon }(\rho _{k},T_{k}^{\ast })}{\nu T_{k,\varepsilon
}^{\ast }\sqrt{T_{k,\varepsilon }^{\ast }}(\rho _{k}D_{k}(T_{k}^{\ast
})-G_{k}(T_{k}^{\ast }))}<\infty ,  \notag
\end{eqnarray}%
independent of $\varepsilon ,$ since $T_{k}^{\ast },T_{k,\varepsilon }^{\ast
}>0$ and $\frac{1}{T_{k,\varepsilon }^{\ast }}\rightarrow \frac{1}{%
T_{k}^{\ast }}.$

Going back to (\ref{94}), we deduce that%
\begin{equation}
\int_{0}^{1}\left\Vert p_{k,\varepsilon }(t)\right\Vert
_{H^{4}(0,L)}^{2}dt\leq \left( \sqrt{2}(1+k^{2})\nu T_{k,\varepsilon }^{\ast
}\left\Vert \overline{p_{k,\varepsilon }^{\prime \prime \prime }}(\cdot
,L)\right\Vert _{V_{1}^{\ast }}\right) ^{2}\int_{0}^{1}\gamma ^{2}(t)dt\leq
C,  \label{98}
\end{equation}%
independent on $\varepsilon ,$ since $\gamma \in L^{2}(0,1).$ Here, $C$
denotes several constants depending on the parameters and $k.$

Therefore, on a subsequence we have 
\begin{equation}
p_{k,\varepsilon }\rightarrow p_{k}\mbox{ weakly in }L^{2}(0,1;H^{4}(0,L)),%
\mbox{ as }\varepsilon \rightarrow 0.  \label{99-1}
\end{equation}%
By the adjoint equation 
\begin{equation*}
E_{0k}(p_{k,\varepsilon })_{t}=T_{k,\varepsilon }^{\ast }F_{0k}^{\ast
}p_{k,\varepsilon }\rightarrow T_{k}^{\ast }F_{0k}^{\ast }p_{k}\mbox{ weakly
in }L^{2}(0,1;H),
\end{equation*}%
hence 
\begin{equation}
(p_{k,\varepsilon })_{t}\rightarrow T_{k}^{\ast }E_{0k}^{-1}F_{0k}^{\ast
}p_{k}:=(p_{k})_{t}\mbox{ weakly in }L^{2}(0,1;H^{2}(0,L)\cap
H_{0}^{1}(0,L)).  \label{99-2}
\end{equation}%
Thus, we get at limit the adjoint equation 
\begin{equation}
-E_{0k}(p_{k})_{t}(t)+T_{k}^{\ast }F_{0k}^{\ast }p_{k}(t)=0\mbox{ a.e. }t\in
(0,1),  \label{99-6}
\end{equation}%
which is (\ref{103}).

Since $H^{4}(0,1)$ is compact in $H^{4-\varepsilon ^{\prime }}(0,1)$ for any 
$\varepsilon ^{\prime }>0,$ it follows by (\ref{99-1}), (\ref{99-2}) and the
Aubin-Lions lemma that 
\begin{equation}
p_{k,\varepsilon }\rightarrow p_{k}\mbox{ strongly in }L^{\infty
}(0,1;H^{4-\varepsilon ^{\prime }}(0,L)).  \label{99-3}
\end{equation}%
By Ascoli-Arzel\`{a} theorem we have 
\begin{equation*}
p_{k,\varepsilon }\rightarrow p_{k}\mbox{ strongly in }C([0,1];H^{2}(0,L)).
\end{equation*}%
Also, (\ref{99-3}) implies 
\begin{equation*}
p_{k,\varepsilon }^{\prime \prime \prime }\rightarrow p_{k}^{\prime \prime
\prime }\mbox{ strongly in }L^{2}(0,1;H^{1-\varepsilon ^{\prime }}(0,L))%
\mbox{ for all }\varepsilon ^{\prime }>0
\end{equation*}%
and also, by the trace convergence 
\begin{equation*}
p_{k,\varepsilon }^{\prime \prime \prime }(\cdot ,L)\rightarrow
p_{k}^{\prime \prime \prime }(\cdot ,L)\mbox{ strongly in }L^{2}(0,1).
\end{equation*}%
On the other hand, by (\ref{97}) 
\begin{equation}
p_{k,\varepsilon }^{\prime \prime \prime }(\cdot ,L)\rightarrow \chi \mbox{
weakly in }V_{1}^{\ast },\mbox{ }  \label{99-4}
\end{equation}%
which combined with the previous convergence yields 
\begin{equation}
\chi (t)=p_{k}^{\prime \prime \prime }(t,L)\mbox{ a.e. }t\in (0,1).
\label{99-5}
\end{equation}%
The convergence of $(v_{k,\varepsilon }^{\ast })_{\varepsilon }$ is ensured
by Theorem 2.2, namely 
\begin{equation*}
v_{k,\varepsilon }^{\ast }\rightarrow v_{k}^{\ast }\mbox{ strongly in }%
C([0,1];H^{2}(0,L))\cap W^{1,2}(0,1;H^{2}(0,L))\cap L^{2}(0,1;H^{4}(0,L))%
\mbox{.}
\end{equation*}

We recall that $w_{k,\varepsilon }^{\ast }\rightarrow w_{k}^{\ast }$ weakly
in $H^{1}(0,1)$ and uniformly in $C([0,1]).$

We have already viewed that the case $\left\Vert w_{k,\varepsilon }^{\ast
}\right\Vert _{V_{1}}<\rho _{k}$ is not acceptable.

Thus, we have to discuss only the case $\left\Vert w_{k,\varepsilon }^{\ast
}\right\Vert _{V_{1}}=\left\Vert \ddot{w}_{k,\varepsilon }^{\ast
}\right\Vert _{V_{1}^{\ast }}=\rho _{k}\sqrt{T_{k,\varepsilon }^{\ast }}.$

We assert that there exists $\alpha _{k}^{\ast }>0$ such that $\alpha
_{k,\varepsilon }\rightarrow \alpha _{k}^{\ast }.$ Indeed, from (\ref{eta})
and (\ref{97}) we derive that $(\alpha _{k,\varepsilon })_{\varepsilon }$ is
bounded. By absurd, if it is not, it means that $\ddot{w}_{k,\varepsilon
}^{\ast }\rightarrow 0,$ that is $\ddot{w}_{k}^{\ast }=0,$ which we have
already seen that it is not compatible with our problem.

Since $\left\Vert \ddot{w}_{k,\varepsilon }^{\ast }\right\Vert _{V_{1}^{\ast
}}=\rho _{k}\sqrt{T_{k,\varepsilon }^{\ast }}$, on the one hand, $\ddot{w}%
_{k,\varepsilon }^{\ast }\rightarrow \ddot{w}_{k}^{\ast }$ weakly in $%
V_{1}^{\ast },$ as $\varepsilon \rightarrow 0.$

On the other hand, by (\ref{eta}) we have that 
\begin{equation*}
\alpha _{k,\varepsilon }\ddot{w}_{k,\varepsilon }^{\ast }\rightarrow
T_{k}^{\ast }\nu \overline{p_{k}^{\prime \prime \prime }}(\cdot ,L)\mbox{
weakly in }V_{1}^{\ast }
\end{equation*}%
because $h_{k,\varepsilon }^{\ast }\rightarrow 0$ strongly in $L^{2}(0,1)$
and so in $V_{1}^{\ast },$ by (\ref{52-6}). We immediately infer that (\ref%
{101}) takes place.

Then, by passing to the limit in (\ref{62}), where $\left\Vert \ddot{w}%
_{k,\varepsilon }^{\ast }\right\Vert _{V_{1}^{\ast }}=\rho _{k}\sqrt{%
T_{k,\varepsilon }^{\ast }}$ we obtain (\ref{102}). This ends the
proof.\hfill $\square $

\section{Problem $(P)$}

\setcounter{equation}{0}

\noindent \textbf{Definition 5.1. }We call a \textit{quasi minimal} solution
to problem $(P)$ a pair $(T^{\ast },w^{\ast })$ given by 
\begin{equation}
T^{\ast }:=\sup\limits_{k\in \mathbb{Z},\mbox{ }k\neq 0}\{T_{k}^{\ast };%
\mbox{ }v_{k}^{\ast }(T_{k}^{\ast })=0\},\mbox{ }w^{\ast
}(t,x)=\sum\limits_{k\in \mathbb{Z}_{k},\mbox{ }k\neq 0}w_{k}^{\ast
}(t)e^{ikx}  \label{T*}
\end{equation}%
where $(T_{k}^{\ast },w_{k}^{\ast })$\ is optimal in $(P_{k}),$\ for each $%
k\in \mathbb{Z}\backslash \{0\}$.

\medskip

We recall that an admissible pair for $(P)$ is a pair $(T_{\ast },w_{\ast })$
with $w_{\ast }\in H^{1}(0,T;L^{2}(0,2\pi )),$ $w_{\ast }(0,y)=0,$ $%
\int_{0}^{T}\int_{0}^{2\pi }\left\vert w_{\ast t}(t,x)\right\vert
^{2}dxdt\leq \rho ^{2}$ and $u(T,x,y)=0,$ $v(T,x,y)=0$.

Further, we shall resume the notation $"\symbol{94}"$ corresponding to the
functions in $(\widehat{P_{k}})$ by the transformation (\ref{t-cap}).

\medskip

\noindent \textbf{Theorem 5.2.} \textit{Let }$(u_{0},v_{0})\in L^{2}(0,2\pi
;H^{3}(0,L)\cap H_{0}^{1}(0,L))\times L^{2}(0,2\pi ;H^{4}(0,L)\cap
H_{0}^{2}(0,L))$ \textit{and} $(u_{0},v_{0})\neq (0,0).$ \textit{If }$(P)$ 
\textit{has an admissible pair }$(T_{\ast },w_{\ast }),$ \textit{there
exists a quasi minimal solution }$(T^{\ast },w^{\ast })$ \textit{to }$(P)$%
\textit{\ given by }(\ref{T*}), \textit{with the corresponding state} 
\textit{\ }%
\begin{eqnarray}
u^{\ast } &\in &W^{1,2}(0,T^{\ast };L^{2}(0,2\pi ;H^{1}(0,L)))\cap
L^{2}(0,T^{\ast };L^{2}(0,2\pi ;H^{3}(0,L)))  \label{sol-u} \\
v^{\ast } &\in &W^{1,2}(0,T^{\ast };L^{2}(0,2\pi ;H^{2}(0,L)))\cap
L^{2}(0,T^{\ast };L^{2}(0,2\pi ;H^{4}(0,L))).  \notag
\end{eqnarray}%
\textit{Moreover, for the modes for which }$\rho _{k},$ $T_{k}^{\ast }$ 
\textit{and} $k$ \textit{satisfy} 
\begin{eqnarray}
&&C_{\nu ,U,k}(1+k^{2})\left\Vert \gamma \right\Vert _{L^{1}(0,T_{k}^{\ast
})}\left( \sqrt{T_{k}^{\ast }}+1\right) <1,\mbox{ }  \label{94-6} \\
&&\rho _{k}\sqrt{T_{k}^{\ast }}\left( 1-C_{\nu ,U,k}(1+k^{2})\left\Vert
\gamma \right\Vert _{L^{1}(0,T_{k}^{\ast })}(\sqrt{T_{k}^{\ast }}+1)\right) 
\notag \\
&>&C_{\nu ,U,k}(1+k^{2})\left\Vert \gamma \right\Vert _{L^{1}(0,T_{k}^{\ast
})}\left\Vert v_{k0}\right\Vert _{H^{4}(0,L)\cap H_{0}^{2}(0,L)},  \notag
\end{eqnarray}%
\textit{the optimal pairs }$(T_{k}^{\ast },w_{k}^{\ast })$ \textit{in }$%
(P_{k})$ \textit{are given by } \textit{\ }%
\begin{equation}
\alpha _{k}^{\ast }T_{k}^{\ast }\ddot{w}_{k}^{\ast }(t)=\nu \overline{%
p_{k}^{\prime \prime \prime }}(t,L),\mbox{ \textit{a.e. }}t\in (0,1),\mbox{ }%
w_{k}^{\ast }(0)=\dot{w}_{k}^{\ast }(T_{k}^{\ast })=0,\mbox{ }\left\Vert
w_{k}^{\ast }\right\Vert _{V_{T_{k}^{\ast }}}=\rho _{k},  \label{94-7}
\end{equation}%
\begin{equation}
\alpha _{k}^{\ast }\rho _{k}^{2}T_{k}^{\ast }+{\rm Re}\int_{0}^{T_{k}^{\ast
}}(\overline{v}_{k}^{\ast }(t),F_{0k}^{\ast }p_{k}(t))_{H}dt=1,\mbox{ }
\label{94-8}
\end{equation}%
\textit{where }$\alpha _{k}^{\ast }>0,$ $v_{k}^{\ast }$ \textit{is the
solution to} (\ref{4})-(\ref{4-3}) \textit{and} $p_{k}$ \textit{is the
solution to the adjoint equation}%
\begin{equation}
-E_{0k}(p_{k})_{t}(t)+F_{0k}^{\ast }p_{k}(t)=0,\mbox{ \textit{a.e.} }t\in
(0,T_{k}^{\ast }).  \label{94-9}
\end{equation}

\medskip

\noindent \textbf{Proof. }By the hypothesis, there exists $T_{\ast }$ such
that $u(T_{\ast },x,y)=v(T_{\ast },x,y)=0.$ Then, by the Parseval identity,
it follows for the modes $k$ that $w_{\ast k}\in H^{1}(0,T_{\ast })$ and $%
u_{k}(T_{\ast },y)=v_{k}(T_{\ast },y)=0.$ Moreover, one can choose $\rho
_{k}:=\left( \int_{0}^{T}\left\vert (w_{\ast k})_{t}(t)\right\vert
^{2}dt\right) ^{1/2},$ which ensures that$\sum\limits_{k\in \mathbb{Z},\mbox{
}k\neq 0}\rho _{k}^{2}\leq \rho ^{2}.$ Thus, $(T_{\ast },w_{\ast k})$ turns
out to be an admissible pair for $(P_{k})$ and correspondingly, $(T_{\ast },%
\widehat{w_{\ast k}})$ is an admissible pair in $(\widehat{P_{k}}).$ Then,
by Theorem 2.3, it follows that $(\widehat{P_{k}})$ has a solution $%
(T_{k}^{\ast },\widehat{w_{k}^{\ast }}),$ and since $T_{k}^{\ast }$ is
minimal it follows that $T_{k}^{\ast }<T_{\ast }.$ The state $\widehat{%
v_{k}^{\ast }}$ has the properties (\ref{vk}) and $\widehat{v_{k}^{\ast }}%
(1,y)=0$. Moreover, by extending $\widehat{w_{k}^{\ast }}$ by $0$ after $%
\widehat{t}=1$ we get $\widehat{v_{k}^{\ast }}(t)=0$ for $\widehat{t}>1$.

Now, we resume the transformation (\ref{t-cap}) and set $\widehat{t}=\frac{t%
}{T_{k}^{\ast }},$ where $t\in \lbrack 0,T_{k}^{\ast }]$ if $\widehat{t}\in
\lbrack 0,1].$ All the results obtained for $(\widehat{P_{k}})$ will be
correspondingly transported to $(P_{k})$ on $(0,T_{k}^{\ast }).$ Thus, $%
(P_{k})$ has an optimal pair $(T_{k}^{\ast },w_{k}^{\ast }),$ $v_{k}^{\ast
}(T_{k}^{\ast })=0$ and $v_{k}^{\ast }(t)=0$ for $t>T_{k}^{\ast }.$

By the third equation in (\ref{2}) we have $u_{k}^{\ast }(t,y):=\frac{i}{k}%
(v_{k}^{\ast })^{\prime }(t,y)$ and so $u_{k}^{\ast }(T_{k}^{\ast })=0$ and $%
u_{k}^{\ast }(t)=0$ for $t>T_{k}^{\ast }.$

We set $T^{\ast }=\sup\limits_{k\in \mathbb{Z},\mbox{ }k\neq 0}\{T_{k}^{\ast
};$ $v_{k}^{\ast }(T_{k}^{\ast })=0\}$ and $w^{\ast }$ as in (\ref{T*}). It
follows that $(T^{\ast },w^{\ast })$ is a quasi minimal solution to $(P)$
and since $T_{k}^{\ast }<T_{\ast }$ it follows that $T^{\ast }<T_{\ast }.$
By the Parseval identity, $u^{\ast }$ and $v^{\ast }$ constructed by (\ref%
{0-3}) satisfy $u^{\ast }(T^{\ast })=v(T^{\ast })=0,$ as required in $(P),$
and $u^{\ast }(t)=v^{\ast }(t)=0$ for $t>T^{\ast }.$

Also, the solution $v_{k}^{\ast }$ to (\ref{4})-(\ref{4-3}) satisfies%
\begin{eqnarray*}
&&\left\Vert v_{k}(t)\right\Vert _{H^{2}(0,L)}^{2}+\int_{0}^{T_{k}^{\ast
}}\left\Vert \frac{dv_{k}}{dt}(t)\right\Vert
_{H^{2}(0,L)}^{2}dt+\int_{0}^{T_{k}^{\ast }}\left\Vert v_{k}(t)\right\Vert
_{H^{4}(0,L)}^{2}dt \\
&\leq &C\left( \left\Vert v_{k0}\right\Vert _{H^{4}(0,L)\cap
H_{0}^{2}(0,L)}^{2}+\int_{0}^{T_{k}^{\ast }}\left\vert w_{k}(t)\right\vert
^{2}dt+\int_{0}^{T_{k}^{\ast }}\left\vert \dot{w}_{k}(t)\right\vert
^{2}dt\right) ,\mbox{ for all }t\geq 0.
\end{eqnarray*}%
Then, 
\begin{eqnarray*}
&&\frac{1}{2\pi }\int_{0}^{2\pi }\left\Vert v^{\ast }(t,x)\right\Vert
_{H^{2}(0,L)}^{2}dx=\sum\limits_{k\in \mathbb{Z},\mbox{ }k\neq 0}\left\Vert
v_{k}^{\ast }(t)\right\Vert _{H^{2}(0,L)}^{2} \\
&\leq &C\sum\limits_{k\in \mathbb{Z},\mbox{ }k\neq 0}\left( \left\Vert
v_{k0}\right\Vert _{H^{4}(0,L)\cap H_{0}^{2}(0,L)}^{2}+\rho
_{k}^{2}(T_{k}^{\ast })^{2}+\rho _{k}^{2}T_{k}^{\ast }\right) \\
&\leq &C\left( \frac{1}{2\pi }\int_{0}^{2\pi }\left\Vert v_{0}(x)\right\Vert
_{H^{2}(0,L)}^{2}dx+\rho ^{2}\max \{(T^{\ast })^{2},T^{\ast }\}\right)
:=I_{1},\mbox{ for all }t\in \lbrack 0,T^{\ast }],
\end{eqnarray*}%
since $\sum\limits_{k\in \mathbb{Z},\mbox{ }k\neq 0}\rho _{k}^{2}\leq \rho
^{2}.$ Similarly, we proceed for showing that $v^{\ast }$ belongs to the
other two spaces.

For $u^{\ast }$ we calculate e.g., 
\begin{eqnarray*}
&&\frac{1}{2\pi }\int_{0}^{2\pi }\left\Vert u^{\ast }(t,x)\right\Vert
_{H^{1}(0,L)}^{2}dx=\sum\limits_{k\in \mathbb{Z},\mbox{ }k\neq 0}\left\Vert
u_{k}^{\ast }(t)\right\Vert _{H^{1}(0,L)}^{2}=\sum\limits_{k\in \mathbb{Z},%
\mbox{ }k\neq 0}\left\Vert \frac{i}{k}(v_{k}^{\ast })^{\prime
}(t)\right\Vert _{H^{1}(0,L)}^{2} \\
&=&\sum\limits_{k\in \mathbb{Z},\mbox{ }k\neq 0}\frac{1}{k^{2}}\left\Vert
(v_{k}^{\ast })^{\prime }(t)\right\Vert _{H^{1}(0,L)}^{2}\leq
\sum\limits_{k\in \mathbb{Z},\mbox{ }k\neq 0}\left\Vert v_{k}^{\ast
}(t)\right\Vert _{H^{2}(0,L)}^{2}\leq I_{1},
\end{eqnarray*}%
and proceed similarly for the other norms, so that $u^{\ast }$ and $v^{\ast
} $ belong to the spaces (\ref{sol-u}).

The second part of the statement follows immediately by making the
transformation $\widehat{t}=tT_{k}^{\ast }$ in (\ref{cond-rok}), (\ref{101}%
)-(\ref{103}). This ends the proof.\hfill $\square $

\medskip

\noindent \textbf{Remark 5.3. }We conclude that the study of the
controllability problem $(P)$ returns the fact that, if $(P)$ has an
admissible pair $(T_{\ast },w_{\ast })$ one can prove that there exists $%
(T^{\ast },w^{\ast })$ which ensures the flow stabilization towards the
stationary laminar regime, with $T^{\ast }\leq T_{\ast }.$ The pair $%
(T^{\ast },w^{\ast })$ has an important property, namely it is constructed
via the solutions of minimal time controllability problems for the modes of
the Fourier transforms of the Navier-Stokes linearized system.

Finally, let us comment the conditions (\ref{94-6}) that enhance the
determination of (\ref{94-7})-(\ref{94-9}) which are the necessary
conditions to be satisfied by the optimal pairs in $(P_{k}).$ One can note
that by (\ref{94-6}) the determination of the optimality conditions can be
done for a finite number of modes $k,$ if $\rho _{k}$, that is, if $\rho $
is large enough.

By strengthening a little bit the first condition in (\ref{94-6}), we have
the possibility of calculating a rough estimate of the number of modes
allowing the determination of the optimality conditions, on the basis of a
relation involving the problem parameters and $T_{\ast }$. Let us assume%
\begin{equation}
C_{\nu ,U,k}(1+k^{2})\left\Vert \gamma \right\Vert _{L^{1}(0,T_{\ast
})}\left( \sqrt{T_{\ast }}+1\right) <1,\mbox{ }  \label{cond-1}
\end{equation}%
which obviously implies the first condition in (\ref{94-6}), since $T_{\ast
}>T_{k}^{\ast }$. However, it provides a smaller number of $k$ than that
obtained by (\ref{94-6}). By calculating the norm $\left\Vert U\right\Vert
_{W^{1,\infty }(0,L)}$ we get $C_{\nu ,U,k}=Ck^{4}\left( \nu +\frac{a}{8\nu }%
(4+L)\right) ,$ with $C$ a positive constant. Thus, we see that a larger
number of $k$ can be obtained function of a smaller controllability cost
reflected by $\left\Vert \gamma \right\Vert _{L^{1}(0,T_{\ast })}$, or a
smaller $T_{\ast }.$ Also, this number depends on $a$ (fixing the
steady-state flow $U(y)),$ $L$ and $\nu ,$ whose values can lead to a
smaller $C_{\nu ,U,k}.$

We also add that one can construct an approximating solution, by restraining
the solution using this finite number of modes. Denote $S_{\rho }:=\{k\in 
\mathbb{Z}\backslash \{0\};$ $k$ satisfies (\ref{cond-1})$\}$ and express
the initial data with respect to these modes only: $u_{0}(x,y)=\sum%
\limits_{k\in S_{\rho }}u_{k0}(y)e^{ikx},$ $v_{0}(x,y)=\sum\limits_{k\in
S_{\rho }}v_{k0}(y)e^{ikx}.$ Then we construct an approximating optimal
state defined by $u^{\ast }(t,x,y)=\sum\limits_{k\in S_{\rho }}u_{k}^{\ast
}(t,y)e^{ikx},$ $v^{\ast }(t,x,y)=\sum\limits_{k\in S_{\rho }}v_{k}^{\ast
}(t,y)e^{ikx},$ corresponding to the optimal pair $(T_{k}^{\ast
},w_{k}^{\ast })$ given by (\ref{94-7})-(\ref{94-9}) with $T^{\ast
}=\sup\limits_{k\in S_{\rho }}\{T_{k}^{\ast };$ $v_{k}^{\ast }(T_{k}^{\ast
})=0\},$ $w^{\ast }(t,x)=\sum\limits_{k\in S_{\rho }}w_{k}^{\ast
}(t)e^{ikx}. $ Such a solution constructed on the basis of a finite number
of modes can be used in numerical computations.

\end{document}